\documentclass[preprint,a4paper,11pt,authoryear]{elsarticle}
\usepackage{endnotes}
\let\footnote=\endnote

\usepackage[round]{natbib}
\usepackage{amsmath}
\usepackage{amssymb}
\usepackage{nicefrac}
\usepackage[font=scriptsize]{subcaption}
\usepackage{adjustbox}
\bibpunct[, ]{(}{)}{,}{a}{}{,}%

\usepackage{comment}

\usepackage{url}
\usepackage{geometry}
\newtheorem{theorem}{Theorem}

\newtheorem{remark}{Remark}

\newtheorem{definition}{Definition}
\newtheorem{assumption}{Assumption}
\newproof{proof}{Proof}

\usepackage{color}
\usepackage[utf8]{inputenc}
\definecolor{mybrown}{rgb}{0.5,0.25,0}
\definecolor{myred}{rgb}{0.75,0,0}
\definecolor{myblue}{rgb}{0,0,0.75}

\newcommand{\ignore}[1]{}

\newcommand{\R}{\mathbb{R}}
\newcommand{\E}{\mathbb{E}}
\journal{European Journal of Operational Research}

\begin{document}

\begin{frontmatter}

\title{Newsvendor Conditional Value-at-Risk Minimisation: a Feature-based Approach under Adaptive Data Selection}

\author{Congzheng Liu\corref{cor1}}
\ead{c.liu19@lancaster.ac.uk}

\author{Wenqi Zhu\corref{cor2}}
\ead{wenqi.zhu@maths.ox.ac.uk}

\cortext[cor1]{Corresponding author}
\address{Department of Management Science,
Lancaster University, Lancaster, LA1 4YX, UK.}

\cortext[cor2]{Co-author}
\address{Mathematical Institute, University of Oxford, Oxford, OX2 6GG, UK.}

\begin{abstract}
The classical risk-neutral newsvendor problem is to decide the order quantity that maximises the expected profit. Some recent works have proposed
an alternative model, in which the goal is to minimise the
\emph{conditional value-at-risk} (CVaR),
a different but very much important risk measure in financial risk management. In this paper, we propose a feature-based non-parametric approach to Newsvendor CVaR minimisation under adaptive data selection (NPC). The NPC method is simple and general. It can handle minimisation with both linear and nonlinear profits, and requires no prior knowledge of the demand distribution. Our main contribution is two-fold. Firstly, NPC uses a feature-based approach. The estimated parameters of NPC can be easily applied to prescriptive analytic to provide additional operational insights. Secondly, unlike common non-parametric methods, our NPC method uses an adaptive data selection criterion and requires only a small proportion of data (only data from two tails), significantly reducing the computational effort. Results from both numerical and real-life experiments confirm that NPC is robust with regard to difficult and large data structures. Using fewer data points, the computed order quantities from NPC lead to equal or less downside loss in extreme cases than competing methods. 
\end{abstract}


\end{frontmatter}


\section{Introduction} \label{se:intro}

In this paper, we focus on \emph{Newsvendor Problem} (NVP), which is a single-period inventory control problem with stochastic demand. In early works on NVP \citep{AHM51,MK51}, it is assumed that the demand in each time period comes from a known probability distribution, and the objective is to
determine the order quantity that maximises the expected profit.

Recently, several works have considered a variant of the NVP in which the
objective is to minimise the \emph{conditional value-at-risk} (CVaR).
The motivation for this is that CVaR is a very much important measure in financial risk management, as pointed out in \cite{RU02}. Compared to traditional NVP, CVaR minimisation NVP offers decision-makers an additional perspective on the worst-case scenario, and it provides decision-makers with a more realistic and nuanced understanding of the risks associated with
inventory management, helping them make more informed and
better decisions. Additionally, the CVaR measure is more appropriate for NVP than other risk measures, such as VaR and TVaR (\emph{tail value-at-risk}), because it is more sensitive to the shape of the tail of the loss distribution. We contend that the CVaR minimisation NVP could be especially advantageous for inventory planners who are risk-averse and seek hedging opportunities using financial instruments, such as call/put options. In the work of \cite{GT07}, a closed-form solution was given for the CVaR-minimisation of NVP. Moreover, a mean-CVaR criterion was considered. Then, \cite{JK07} proposed an extended model where the inventory manager can control internal and customer-oriented performance measures. \cite{CXZ09}, later on, investigated the optimal pricing and ordering decisions in a single framework. This risk averseness in pricing competition is further studied in \cite{WZT14}. Other relevant literature can be found in \cite{AS17,XZRMP15,WZT13,CWW09}. 

The CVaR minimisation NVP is one form of the risk-averse NVP. Therefore, it inherits the features of the risk-averse NVP. The versatility of the risk-averse NVP family has also been extensively studied \citep[see][for example]{WWS09,CR08,CR11}. In real-world scenarios, there are numerous instances where stakeholder decisions do not align with the expected profit maximisation order quantity \citep{SC00}. Thus, the assumption of risk neutrality is not always applicable. Moreover, the expected profit foregone by utilising a risk-averse quantity versus the conventional risk-neutral quantity may vary depending on the decision-maker's risk tolerance level. Hence, comparing the two quantities, and corresponding expected profits, may not be practical as they apply different objectives from the outset.

In all of the above-mentioned works, it is assumed that the demand comes from a known family of probability distributions with correct parameters.
In real life, unfortunately, model correctness is rarely assured. Assuming that historical data is available, one can attempt to address this issue by decomposing the problem into a forecasting phase and an optimisation phase, commonly called the disjoint approach \citep{LLS22} or SEO approach \citep{BR19}. However, if the forecasting model is misspecified, and/or there is substantial noise in the data, then these factors may impact the optimisation phase in an unexpected way. This could possibly lead to sub-optimal solutions or even nonsensical negative order decisions. Moreover, this issue becomes even more severe in the CVaR minimisation of NVP.  Given that the CVaR concerns observations with extreme values, which are often treated as outliers in traditional statistical approaches, the computed order quantities could underestimate the downside risk and lead to significant losses in extreme cases \citep{GSU03,YWL13}. The forecasting accuracy can be slightly improved by considering feature-based demand data \citep{Va98a}, or by using an alternative statistical approach, such as bootstrapping \citep{ET94} or extreme value theory \citep{DFF06}. Yet, the performance depends heavily on the form of the profit function. 

To overcome these difficulties, we propose a \emph{feature-based non-parametric approach of CVaR minimisation} (i.e.,``NPC" for short). Our proposed approach does not make assumptions about the demand distribution, but works directly with historical data and considers features related to the demand. 
As such, the NPC is very robust with regards to the data structures and is adaptive to different forms of the profit function, including both linear and nonlinear. Moreover, our NPC model can cope with different forms of demand-feature relationships. Additionally, the estimated parameters of NPC can be explained statistically and readily applied to prescriptive analytic to provide additional operational insights. Note that there are other non-parametric approaches for classical risk-neutral NVP including SAA \citep{LPU15}, SPO loss \citep{EG17}, NV-features \citep{BR19} and IMEO \citep{LLS22}. However,  to our best knowledge, there are no literature in applying the non-parametric approaches to NVP under CVaR minimisation. 

Another key feature of the NPC model is that we use an adaptive sampling criterion to select the data points for the CVaR risk minimisation process. Unlike common non-parametric approaches, which become extremely computationally expensive under large instances, our proposed NPC approach is ``smarter" and ``faster" as it only requires a small proportion of data (only data from two tails) for the minimisation.  We also give a rigorous proof that under suitable assumptions, the minimisation result from the reduced data set model well approximates the result from the model with full data set.  Through extensive experiments, on both hypothetical and real data, we verify that with a small proportion of well-selected data points, the computed order quantities from our NPC method lead to equal or less downside loss in extreme cases than competing methods. 

The paper is organised as follows. We review some well-known results of the single-period NVP in Section \ref{se:news}. Then, in Section \ref{se:feature} and Section \ref{se:adaptive}, we present the method of NPC in detail. In those two sections, we discuss the two distinct features of our NPC model, namely feature-based estimation and adaptive data selection, respectively. We also provide theoretical analysis to support the adaptive data selection. Sections
\ref{se:art} and \ref{se:ex_exp} give computational results on hypothetical data, with baseline experiments and extended experiments. Section \ref{se:exp_real}, on the other hand, applies NPC to a real-life example within a food bank. Finally, Section \ref{se:conclusion} contains some
concluding remarks.


\section{Reviews on Newsvendor Problem} \label{se:news}

In this section, we review the Newsvendor problem. In subsection \ref{sub:class}, the classic NVP that maximises the expected profit is considered. In subsection \ref{sub:nnvp}, we define the profit function of a nonlinear NVP. Subsection \ref{sub:CVaR} gives a closed-form solution to the NVP under CVaR minimisation. 

\subsection{Classic newsvendor} \label{sub:class}

In the simplest NVP, a company
purchases goods at the beginning of a time period, and aims to sell them by the end of the period. The demand during the period is a random variable $\tilde{d}$ with a known distribution function $F$. We also know that $c$ is the cost of purchasing one unit of item, $r$ is the revenue gained by selling one unit of item, $v$ is the disposal cost of each unsold unit of item and $g$ is the shortage cost of each unit of unsatisfied demand for item. We note that salvage sales and emergency backup can lead to negative values for $v$ and $g$, respectively.
We assume without loss of generality that $r > c \ge 0$ and $c > - \, v$. The retailer must decide how many units of the item to order before the
start of the sales period. We let $x$ denote the number of units
ordered and we assume
for simplicity that $x$ is continuous. For a given value of $x$, and a given realisation $d$
of $\tilde{d}$, the profit over the period is:
\begin{equation} \label{eq:nvp}
\pi(x,d) := r\,\min\{x,d\} - c\,x -v\,[x-d]^{+} -g\,[d-x]^{+}.
\end{equation}
In the classical NVP, the goal is to find a value for $x$ that maximises the total expected profit, which can be given in closed form \citep{AHM51,Ch12}:
\begin{equation} \label{eq:nvpsol}
x^* = F^{-1}\left(\frac{U}{E+U}\right),
\end{equation}
where $F^{-1}$ is the inverse of the distribution function $F$, $E := c+v$ denotes the overage cost, and $U := r-c+g$ denotes the underage cost. 

\subsection{Nonlinear newsvendor} \label{sub:nnvp}

In the general nonlinear NVP, as defined in \cite{LLS22}, the profit function takes the form:
\begin{equation} \label{eq:nnvp}
    \pi(x,d) :=
    \begin{cases}
        R(x,d)-C(x,d)-V(x,d),& \text{for } x \geq d\\
        R(x,d)-C(x,d)-G(x,d),& \text{for } x< d,
    \end{cases}
\end{equation}
where $R$, $C$, $V$ and $G$ are now \emph{functions} rather than constants.

The nonlinear NVP can be seen as an extension of the classical NVP, as it enables one to model more real-life problems, e.g., with nonlinear shortage cost due to the damage of reputation. The detailed motivation can be found in \citet{PK91,Kh95,LLS22}. In general, however, a closed-form expression as a quantile is unlikely to exist for nonlinear NVP. In such cases, one has to resort to numerical integration or simulation techniques to solve the problem. Due to this reason, we note that the current existing parametric approach is not directly applicable to the CVaR minimisation of nonlinear NVP.

\subsection{NVP under CVaR minimisation} \label{sub:CVaR}

\begin{figure}[htb]
    \centering
    \caption{The cumulative distribution function of $L(x,\tilde{d})$ and the $\beta$-tail distribution. In this plot, $\beta = 95\%$}.
    \label{fig:cdf}
    \begin{subfigure}[b]{0.47\textwidth}
    \centering
    \includegraphics[width=0.8\textwidth]{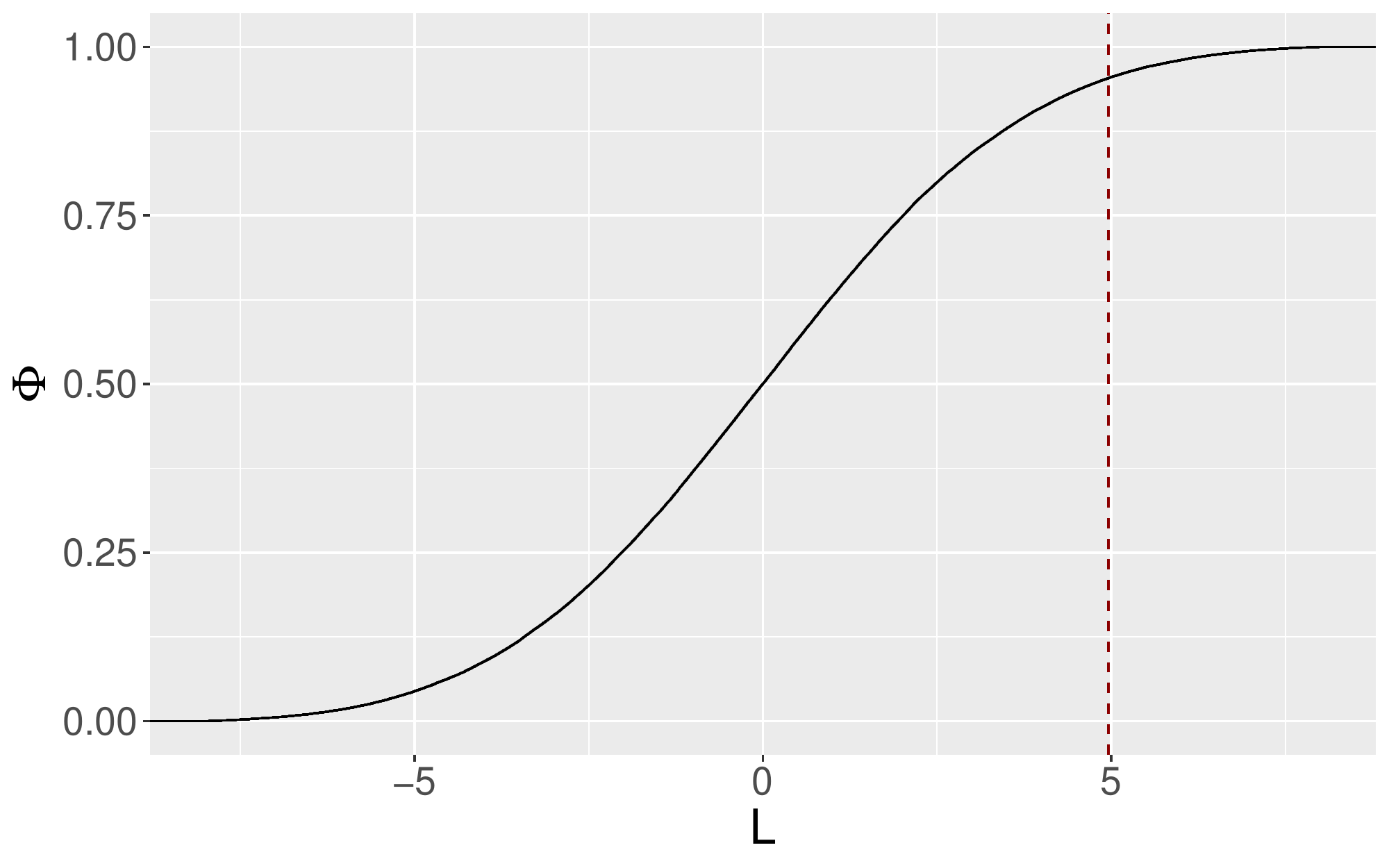}
    \caption{cumulative distribution}
    \end{subfigure}
    \hfill
    \begin{subfigure}[b]{0.47\textwidth}
    \centering
    \includegraphics[width=0.8\textwidth]{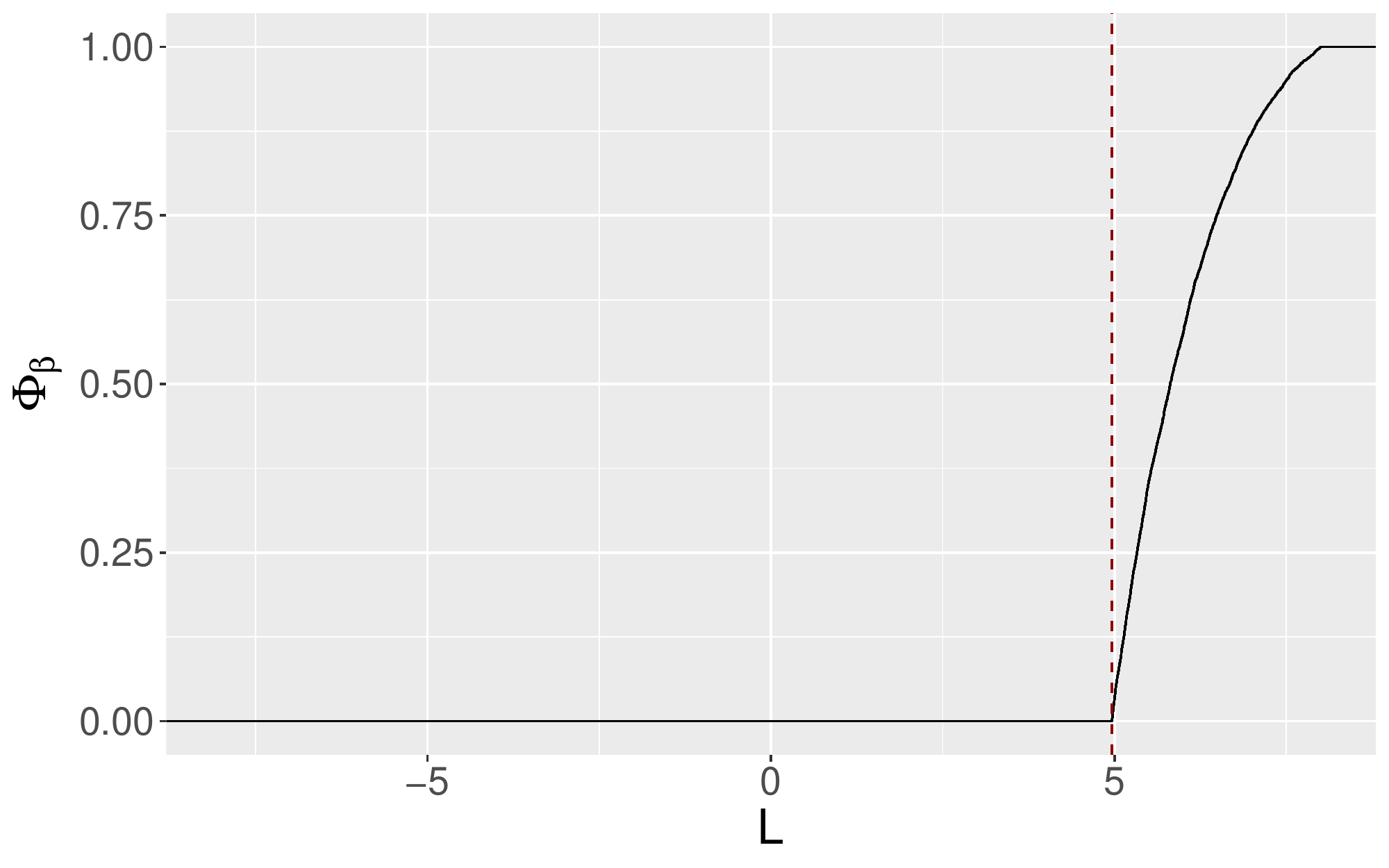}
    \caption{$\beta$-tail distribution}
    \end{subfigure}
\end{figure}

We define the loss function as 
\begin{equation}
    L(x,\tilde{d}) := - \pi (x,\tilde{d})
\end{equation}
and the magnitude of the loss for a given realisation $d$ of $\tilde{d}$ and a fixed $x$ is represented by $L(x,d)$. Let $\Phi (\eta|x):= \mathbb{P} \{L(x,\tilde{d}) \leq \eta\}$ denote the cumulative distribution function of $L$. For $ \beta \in [0, 1)$, we define the $\beta$-VaR of the distribution by 
\begin{equation} \label{eq: alpha}
    \alpha(x, \beta):= \inf_{\eta \in \R}\{\eta|\Phi (\eta|x) \geq \beta\}.
\end{equation} 
In other words, $\alpha$ is the risk threshold such that the probability of incurring a loss greater than $\alpha$ is exactly $(1-\beta)$. If the cumulative distribution function is continuous and the infimum is reached, we can write  $P \{L(x,\tilde{d}) \ge \alpha\} = 1-\beta$ and $\Phi(\alpha | x ) = P \{L(x,\tilde{d}) \le \alpha\} = \beta$.
A $\beta$-tail distribution function that focuses on the upper tail part of the loss distribution can be formed as \citep{RU02}:
\begin{equation}
\Phi_\beta (\eta|x) := 
    \frac{\Phi(\eta|x)-\beta}{1-\beta}, \quad \text{for } \eta \geq \alpha(x, \beta).
\end{equation}
We plot an illustrative distribution function of $\Phi(\eta|x)$ and $\Phi_{\beta}(\eta|x)$ in Figure \ref{fig:cdf}. It is easy to see that the $\beta$-tail distribution is formed by picking the top $(1-\beta)$ proportion of $\Phi(\eta|x)$ values, scaling those values by an affine transformation and setting the rest of the $\Phi(\eta|x)$ values to $0$.

Using the $\beta$-tail distribution, \citet{RU02} formulated the \emph{$\beta$-conditional value-at-risk ($\beta$-CVaR) risk minimisation model} as, 
\begin{equation} \label{eq:aufunction}
\min_{x, \alpha} F_\beta (x,\alpha), \qquad F_\beta (x,\alpha):= \alpha + \frac{1}{1-\beta} \E \left[[L(x,\tilde{d})-\alpha]^+ \right]. 
\end{equation}
This minimisation problem is also the key model that we will focus on in this paper. For details and extensions see \cite{RU02,Mi14,HNS15}.

The closed-form solution to the CVaR minimisation of a (linear) NVP is derived by \citet{GT07},
\begin{equation}
\label{eq:CVaR}
    x^* = \frac{E+W}{E+U}F^{-1}\left(\frac{U(1-\beta)}{E+U}\right) + \frac{U-W}{E+U}F^{-1}\left(\frac{E\beta +U}{E+U}\right),
\end{equation}
where $E := c+v$, $U := r-c+g$, and $W := r-c = U-g$. In particular, when $g = 0$, we have a simpler result:
\begin{equation} \label{eq:simCVaR}
    x^* = F^{-1}\left(\frac{U(1-\beta)}{E+U}\right).
\end{equation}

\begin{remark}
It is worth noticing that (\ref{eq:nvpsol}), (\ref{eq:CVaR}) or (\ref{eq:simCVaR}) depends only on two parameters $g$ and $\beta$. In particular, when $g=0$, the difference between  (\ref{eq:nvpsol}) or (\ref{eq:simCVaR}) is only the coefficient in the argument of the inverse $F^{-1}$. Moreover, when $\beta=0$, (\ref{eq:CVaR}) reduces to the classical expected profit maximisation solution in (\ref{eq:nvpsol}). This result shows the definition of the $\beta$-CVaR is consistent. More details can be found in  Figure \ref{fig:mean_CVaR} in \ref{app:e_CVaR}.
\end{remark}


\section{Feature-based estimation of NPC} \label{se:feature}

Although the closed form solution for linear NVP under CVaR minimisation is available, in practice, the demand distribution $\tilde{d}$ and the underlying distribution of $L(x, \tilde{d})$ are often unknown \emph{a priori}, not to mention the difficulty of solving nonlinear NVP. As such, the closed-form solutions in \eqref{eq:CVaR} and \eqref{eq:simCVaR} are usually not directly applicable.
Instead of computing the closed-form solution analytically, one could perform a fine approximation using non-parametric approaches. 

A straightforward, and commonly applied, non-parametric approach is to use the sample average approximation (SAA) to form an empirical feature-based risk \citep{LPU15}. That is;  given historical demand observations $[d_1, \dots, d_s]$, the SAA aims to find an order quantity $x$ that
\begin{equation} \label{eq:datanvp}
\min_{x, \alpha} \left(\alpha + \sum_{t =1}^s \frac{[L(x, d_t)-\alpha]^+}{ (1-\beta)s } \right).
\end{equation}

However, we note that this method only outputs a point order decision. Moreover, the minimisation in \eqref{eq:datanvp} cannot compute accurate inventory information that is related to a strong trend or seasonality of the demand,  both of which are highly important in many real-world inventory
control problems.

To overcome these drawbacks, in this paper, we design a \emph{feature-based model} to estimate the empirical risk instead. We collect data on features of the demand as well as the demand itself. We assume that the historical data are
$[(\mathbf{z}_1,d_1),\dots,(\mathbf{z}_s,d_s)]$. For $t=1, \dotsc, s$, each $\mathbf{z}_t := [z^1_t,\dots,z^p_t]$ represents $p$ features related to the
demand, including both exogenous information (prices, promotions) and intrinsic information (seasonal patterns, lagged demand). We consider
$x = h (\mathbf{z})$. Consequently, the CVaR minimisation problem becomes
\begin{equation} \label{eq:ex_empirical}
 \min_{h, \alpha} \tilde{F}_{\beta}(h, \alpha), \qquad \tilde{F}_{\beta}(h, \alpha): =\alpha + \frac{1}{(1-\beta)s} \sum_{t =1}^s {[L(h(\mathbf{z}_t), d_t)-\alpha]^+}.
\end{equation}
We refer to  $ \tilde{F}_{\beta}(h, \alpha)$ as \emph{the empirical feature-based risk}. Both $h$ and $\alpha$ are optimised in the CVaR minimisation. We then evaluate $x$ in the proceeding period, using $x_{s+1} = h(\mathbf{z}_{s+1})$. In the simplest case, $h$ can be of the linear form:
\begin{equation}\label{function h}
    h(\mathbf{z}_t) := \mathbf{z}_t^{\mathsf{T}}\boldsymbol{\gamma} = \sum_{j=1}^p z_t^j \gamma^j
\end{equation}
where $\boldsymbol{\gamma} \in \R^p$. In this linear model, the optimised $\boldsymbol{\gamma}$ can be interpreted as the `effective ratio' of given features. Namely, the order should be increased by $\gamma^j$ units to achieve a minimum CVaR if the feature $z^j$ is increased by one unit. Using a feature-based model, not only can we find the optimal order that maximises the profit or minimises the CVaR, but we can also extract important information (such as price or seasonality) that affects the profit. We give illustrative examples in Sections \ref{se:art} and \ref{se:ex_exp}, where the results of NPC method and traditional methods are compared. Then, in Section \ref{se:exp_real}, we show how insights can be drawn from the NPC results. Moreover, for more complicated demand structures, under suitable assumptions, we can use a nonlinear representation of $h$ in the minimization (for instance, higher-order polynomial representation, or log-linear representation). These assumptions are discussed in Remark \ref{remark: valid function space} in Section \ref{sub:fun_space}.

While additional forms of $h$ can prove useful in certain situations, the experiments conducted in Section \ref{se:art} and Section \ref{se:ex_exp} suggest that the linear form of $h$ is generally adequate for accurate estimation. The idea of feature-based newsvendor stems from the \emph{empirical risk minimisation principle} that has been widely adopted by the machine
learning community. Related works include \cite{EG17,BR19,LLS22}. Yet, it seems that there are no studies conducted on feature-based NVP under CVaR minimisation. Furthermore, we argue that the features crucial to CVaR minimisation may be very different from those of expected profit maximisation.


\section{Adaptive data selection of NPC} \label{se:adaptive}

In this section, we propose an adaptive way of selecting the data for NPC. Instead of minimising empirical feature-based risk using the whole data set $\{(\mathbf{z}_t, d_t)\}_{1 \le t \le s}$, we carefully select a $2\times(1-\beta)$ portion of the data and use the reduced data set to minimise a tailored risk function. Note that we require $\beta > 50\%$; the standard value of $\beta$ is $90\%$ or $95\%$ in practice. 

The section is ordered as such; In Subsection \ref{sub:method}, we present the adaptive data selection of NPC in detail. In Subsection \ref{sub:proof} -- \ref{sub:empirical_risk}, we prove that the risk generated by our model converges to the true feature-based risk. We note that $h$ denotes a general choice of function throughout this section.

\subsection{Data selection strategy in NPC}
\label{sub:method}

Assuming that $\{\boldsymbol{z}_t, d_t\}_{1 \le t \le s}$ are samples drawn from a true multi-variant data distribution, and the observed time series $[d_1,\dots,d_s]$ can be decomposed into the systematic component $T$ and the irregular component $\epsilon$. After the decomposition, the set $\{d_t\}_{1\le t \le s}$ corresponds to a set of systematic components $\{T_t\}_{1\le t \le s}$ and a set of irregular components $\{\epsilon_t\}_{1\le t \le s}$. For simplicity, we assume the irregular components are independent and identically distributed (i.i.d.). We can re-write the loss as
\begin{equation}
    L(h(\mathbf{z}_t), d_t) = \tilde{L}(h(\mathbf{z}_t), T_t, \epsilon_t)  \qquad \text{ for all } 1\le t \le s 
\end{equation}
We observe that in the NVP, the loss takes a large value if and only if the irregular component $\epsilon_t$ takes extreme absolute values. This observation motivates us to design an adaptive selection criterion. 

We define the `worst' scenarios as the `smallest' $(1-\beta)$ proportion and the `largest' $(1-\beta)$ proportion of the data in regard to their irregular components $\epsilon$. We denote the selected irregular components
in ascending order as
\begin{equation}
\mathcal{E} := \{\epsilon_{i_1}, \dotsc, \epsilon_{i_m}, \epsilon_{i_{m+1}}, \dotsc,\epsilon_{i_{2m}}\}
\label{data select 1}
\end{equation}
where $m : =\lceil (1-\beta)s \rceil$. The first $m$ items of  $\mathcal{E}$ are the $m$ smallest $\epsilon_t$ values, and the last $m$ items are the $m$
largest. We also denote the index set of the chosen data as
\begin{equation}\label{eq: M}
    M := \{i_1, \dotsc, i_{2m}\}.
\end{equation} 
 By definition, the cardinality of $|M| = 2m = 2\lceil (1-\beta)s \rceil$. 
 
To measure CVaR, We aim to select the $(1-\beta)$ portion of the largest loss function $\tilde{L}(h(\mathbf{z}_t), T_t, \epsilon_t)$. It is important to note that the extreme loss may be caused by either the  `smallest'  $(1-\beta)$ proportion or the `largest'  $(1-\beta)$ proportion of the irregular components in data, depending on the specific loss function.  A numerical illustration is provided in \ref{app:sm}. Therefore, we need to select a total of $2\lceil (1-\beta)s \rceil$ `largest' or  `smallest'  sample points to generate the largest $(1-\beta)$ proportion of the loss function.
Formally, we define the (discrete) risk threshold in Definition \ref{def 1}.
\begin{definition}
Let $\alpha$ be the (discrete) \textbf{risk threshold} such that there exist exactly $\lceil (1-\beta)s \rceil$ values of the loss $\tilde{L}$ which have a larger value than $\alpha$. We denote the index set of these data as 
\begin{equation}\label{eq: S}
    S := \{j_1, \dotsc, j_{m}\}.
\end{equation} 
The cardinality of $|S| = m = \lceil (1-\beta)s \rceil$. 
\label{def 1}
\end{definition}

Employing the aforementioned definitions and the adaptive data selection technique, we define the tailored risk as follows,
\begin{equation} \label{eq:empirical}
 \min_{h,\alpha} \hat{F}_{\beta}(h, \alpha),\quad \hat{F}_{\beta}(h, \alpha) :=  \alpha + \frac{1}{m}\sum_{t \in M} \big[ \tilde{L}(h(\mathbf{z}_t), T_t, \epsilon_t)-\alpha\big]^+.
\end{equation}
Note that although $\mathcal{E}$ contains $2m$ data points with extreme values in the irregular component, only $m$ of these data points result in a loss function such that $\tilde{L}(h(\mathbf{z}_t), T_t, \epsilon_t) \ge \alpha$. Therefore, only $m$ terms make a positive contribution to the summation in equation \eqref{eq:empirical}. As a result, when calculating the tailored risk, the denominator is $m$.

In the subsequent sections, we provide a comprehensive analysis that our tailored risk, $\hat{F}_{\beta}$, converges to the empirical risk, $\tilde{F}_{\beta}$, which in turn converges to the true risk, $F_{\beta}$. By using only a proportion of the available data, our proposed risk model improves the computational efficiency of the minimisation process. 

\subsection{Proof of convergence}
\label{sub:proof}

In this subsection, we prove that under suitable assumptions, the tailored risk obtained from our carefully selected $2\lceil (1-\beta)s \rceil$ proportion of data approximates the true feature-based risk in an appropriate function space.

We prove the convergence in two parts. 
\begin{itemize}
    \item \textbf{Part I: Tailor risk converges to empirical risk. }We prove that, under suitable assumptions and for each given data set, there is an equivalence relationship between the tailored risk $ \hat{F}_{\beta}(h, \alpha) $ and the empirical feature-based risk  $\tilde{F}_{\beta}(h, \alpha)$. Note that the empirical risk defined in \eqref{eq:ex_empirical} is calculated using $s$ data points, whereas the tailored risk defined in \eqref{eq:empirical} is calculated using $2\lceil (1-\beta)s \rceil$ selectively chosen data points. We highlight that this step gives crucial theoretical support to our selection criteria. It is also one of the most mathematically demanding parts of the paper. To prove this, we require mild assumptions for the loss function.
    
     \item \textbf{Part II: Empirical risk converges to true risk.} 
     We prove that, under a suitable function space $\mathcal{F}$ for the estimator, the empirical risk converges to the true risk, $ \Tilde{F}_{\beta}(h, \alpha) \rightarrow {F}_{\beta}(h_{\mathcal{F}}, \alpha)$, as sample size grows. Note that we define the true feature-based risk as $F_\beta (h_\mathcal{F},\alpha)$, where $h_\mathcal{F}$ is the best estimator function in the function space $\mathcal{F}$. This represents the risk calculated using the best estimator function in $\mathcal{F}$ on the actual distribution of data points. This step is contingent upon the theoretical postulations and scholarly discourse stemming from the field of statistical learning theory. To prove this, we require assumptions on the function space $\mathcal{F}$ of the estimator function $h$. 
\end{itemize}
Using \textbf{Part I} and \textbf{II}, we deduce our main convergence result, $ \hat{F}_{\beta}(h, \alpha) \rightarrow {F}_{\beta}(h_{\mathcal{F}}, \alpha)$, as sample size grows. That is, provided with a large data set, the tailored risk well approximates the true risk in the function space $\mathcal{F}$. 

\subsection{Part I:  Tailor risk converges to empirical risk}
In this part, we prove the validity of the adaptive data selection criterion. The key idea of our proof is twofold. Firstly, we demonstrate that when we compute the empirical feature-based risk, only the data $(T_t, \epsilon_t)$ corresponding to $\tilde{L} \ge \alpha$ contribute to the expectation, while the remaining data points result in a zero term in the expectation. Secondly, using our adaptive data selection criterion, the $2\lceil(1-\beta)s\rceil$ selected indices in $M$ are adequate to cover the data points that generate a non-zero expectation in $\tilde{F}_{\beta}(h, \alpha)$. To complete these proofs, we need the following assumptions. 

\begin{assumption} \label{assumption 1}
Assume the distribution of the irregular component and the loss function satisfies: 
\begin{enumerate}
    \item \textbf{Distribution Assumption for $\epsilon$:} $\{\epsilon_t\}_{1\le t\le s}$ are i.i.d. random variables from a distribution with zero mean. Let $\Phi_{\epsilon}(\eta) := \mathbb{P}(\epsilon<\eta_{\epsilon}) $ be the cumulative distribution function of the distribution, such that 
    $\displaystyle \lim_{\eta_{\epsilon} \to -\infty} \Phi _{\epsilon}(\eta) = 0$ and $\displaystyle \lim_{\eta_{\epsilon} \to \infty} \Phi_{\epsilon}(\eta) = 1$.
    \item  \textbf{Continuity Assumption:} Let $( \mathcal{F} ,\mathcal{T}, (-\infty, \infty))$ be the feasible set for the distribution of $(h,  T, \epsilon)$.  The loss function $\tilde{L}$ is continuous with respect to $\epsilon$ for all $(h,  T) \in (\mathcal{F} ,\mathcal{T})$. 
    \item  \textbf{Tail Assumption:} For all $(h,  T) \in (\mathcal{F},\mathcal{T})$, we assume that one or both tail of the loss function is monotonic and unbounded as $|\epsilon| \rightarrow \infty$. Namely, one of the below tail scenarios is true: 
\begin{eqnarray*}
\lim\limits_{\epsilon\rightarrow -\infty} \tilde{L}(h,  T,\epsilon) \rightarrow \infty \text{monotonically and } \lim\limits_{\epsilon\rightarrow \infty} \tilde{L}(h,  T, \epsilon) \text{ bounded, or,}
\\ \lim\limits_{\epsilon\rightarrow \infty} \tilde{L}(h,  T, \epsilon) \rightarrow \infty \text{ monotonically and } \lim\limits_{\epsilon\rightarrow -\infty} \tilde{L}(h, T, \epsilon) \text{ bounded, or,}  
\\ \lim\limits_{\epsilon\rightarrow \pm\infty} \tilde{L}(h,  T, \epsilon) \rightarrow \infty \text{ monotonically}.
\end{eqnarray*}
\end{enumerate}
\end{assumption}

\begin{theorem}
Let $\alpha$ be the risk threshold as defined in Definition  \ref{def 1}, and consider the adaptive data selection criterion given by Equation \eqref{data select 1} for a given data set. Under Assumption 1, for a suitable $\beta$, we have $\hat{F}_{\beta}(h, \alpha)   = \tilde{F}_{\beta}(h, \alpha)$.
\label{thm2}
\end{theorem}
\proof
In the proof, we denote $\tilde{L}(h, T, \epsilon)$ as $\tilde{L}(\epsilon)$ for notational simplicity.

\textbf{Step 1: Region for extreme values of the loss function.} Under the tail assumption, there exist positive constants $L_a, \epsilon_a \in \R$ satisfying the following conditions:
\begin{equation}
\begin{cases}
    |\tilde{L}(\epsilon)| \ge L_a, & \quad \text{ for } \, |\epsilon| \ge \epsilon_a, \\
    |\tilde{L}(\epsilon)| < L_a, & \quad \text{ for }  \, - \epsilon_a \le \epsilon \le \epsilon_a.
\end{cases}
\label{cond2}
\end{equation}
Also, $\tilde{L}(\epsilon)$ exhibits monotonic behaviour in the regions $(-\infty, -\epsilon_a]$, $[\epsilon_a, \infty)$, or $(-\infty, -\epsilon_a] \cup [\epsilon_a, \infty)$, contingent upon the three distinct tail scenarios.

\textbf{Step 2: Distribution for $\epsilon$.} Consider $\{\boldsymbol{z}_t, T_{t}, \epsilon_{t}\}_{1 \le t \le s}$ consisting of i.i.d. samples drawn from a true multi-variant data distribution. In light of the assumption on the distribution for $\epsilon$ , when the sample size is sufficiently large, there exists a $\beta\in (0, 1)$ such that
\begin{equation*}
\epsilon_{i_1} \le \dots \le \epsilon_{i_m} \le -\epsilon_a \le \epsilon_a \le \epsilon_{i_{m+1}} \dots \le \epsilon_{i_{2m}}
\end{equation*}
where $m = \lceil (1-\beta)s \rceil$, and $M := \{i_1, \dots, i_{2m}\}$ comprises either the `largest' $m$ values or the `smallest' $m$ values of the irregular component.

\textbf{Step 3: Relationship between $S$ and $M$.}
We aim to prove $S \subseteq M$ where $S$ and $M$ are defined in \eqref{eq: S} and \eqref{eq: M}, respectively. We assume for contradiction that there exists an ${i_*} \in S$ and ${i_*} \notin M$. The assumption ${i_*} \notin M$ implies that $\epsilon_{i_*}$ is neither among the `largest' nor `smallest' $m$ values of the irregular component,  which leads to the following inequality:
\begin{eqnarray*}
    \epsilon_{i_1} \le \dotsc \le   \epsilon_{i_m} \le \epsilon_{i_*} \le \epsilon_{i_{m+1}} \dotsc \le  \epsilon_{i_{2m}}.
\end{eqnarray*}
According to \eqref{cond2} and the monotonicity property, depending on the three different tail scenarios, we have:
\begin{eqnarray*}
    &&\tilde{L}( \epsilon_{i_*}) \le  \tilde{L}( \epsilon_{i_m}) \le    \dotsc \le   \tilde{L}(\epsilon_{i_1}),   \text{ or }
     \\    &&\tilde{L}( \epsilon_{i_*}) \le \tilde{L}(\epsilon_{i_{m+1}})  \le    \dotsc \le     \tilde{L}(\epsilon_{i_{2m}}) , \text{ or both}. 
\end{eqnarray*}
In all cases, $\tilde{L}(\epsilon_{i_{m+1}})$ is not one of the largest $m$ values of the loss. However, $S$ is defined as the collection of $m$ largest values of the loss. Consequently, we find that $i_* \notin S$, which contradicts our assumption. Thus, we conclude that $S \subseteq M$. A numerical illustration is provided in \ref{app:sm}.

\textbf{Step 4: Equivalence.} 
Finally, we prove the relationship between  $\hat{F}_{\beta}(h, \alpha)$ and $\tilde{F}_{\beta}(h, \alpha)$. Using the relationship $S \subseteq M$, we can partition $ M$ as $M = S \cup (M/S)$, it follows from definition that
\begin{eqnarray*}
 &&\hat{F}_{\beta}(h, \alpha)  \underset{\eqref{eq:empirical}}{=}   \alpha+ \frac{1}{m}\sum_{t \in M}\big[\tilde{L}(\epsilon_t) -\alpha)\big]^+  
 \\ && \quad = \alpha+ \frac{1}{m} \bigg[\sum_{t \in S}\big[\tilde{L}(\epsilon_t)-\alpha)\big]^+  + \underbrace{ \sum_{t \in M/S}\big[\tilde{L}(h(\epsilon_t)-\alpha)\big]^+}_{=0\text{ by Definition \ref{def 1}} }\bigg]
 \\ && \quad = \alpha+ \frac{1}{m}\sum_{t \in S}\big(\tilde{L}(h(\epsilon_t) -\alpha\big) \underset{\eqref{eq:ex_empirical}}{=} \tilde{F}_{\beta}(h, \alpha).
\end{eqnarray*} 
\endproof
\vspace{-0.3cm}
\hfill $\square$

The proof of Theorem \ref{thm2} sheds light on the significance of the assumptions made for the loss function. It also confirms the universality of our adaptive data selection method for CVaR problems and other risk management scenarios. To elaborate, the assumption of monotonicity for the entire domain of the loss function is not required since we are only dealing with $(1-\beta)$ proportion of extreme values in the loss function and $2(1-\beta)$ proportion of extreme values in $\epsilon$. Instead, we only need the loss function to have monotonic tails. As demonstrated in Step $1$ of the proof, this is sufficient to establish \eqref{cond2} and ensure the validity of the proof. Such properties enhance our flexibility in selecting loss functions and broaden the range of minimisation problems to which our method can be applied.

Step $3$, which involves proving that $S \subseteq M$, is crucial for understanding the methodology of the proof. Step $3$ suggests that if the loss function is heavily influenced by extreme inputs (i.e., $|L| \rightarrow \infty$ as $|\epsilon| \rightarrow \infty$), then selecting $2\lceil (1-\beta)s \rceil$ data points from $M$ is sufficient to generate the cases with the largest $(1-\beta)$ proportion of the loss.

In general, whenever the objective is to minimise the risk of extreme losses and the loss function satisfies the tail assumption, our adaptive data selection of NPC can be considered an efficient data selection method for improving computational efficiency.

\subsection{Part II:  Empirical risk converges to true risk}
\label{sub:empirical_risk}

Now, we show the empirical risk converges to the true risk. In other words, we would like to investigate the behaviour of the estimator $h$ as sample sizes $s$ increase. Therefore, in this section, we use $h^{(s)}$ to denote the  estimator constructed by $s$ data points. 

\subsubsection{For a fixed $h$}

First, we show that the empirical risk converges to the true risk (expected risk over the entire data distribution) for a fixed $h \in \mathcal{F}$ as the sample size grows. We give the formal definition of consistency following \cite{VS11}. 

\begin{definition}
The data points $\{\boldsymbol{z}_t, d_t\}_{1 \le t \le s}$ are i.i.d. samples drawn from a true multi-variant data distribution. For each $s \in \mathbb{N}$, we denote
$h^{(s)}$ be an estimator constructed by the first $s$ data points. The estimation is called consistent with the true multi-variant data distribution in $\mathcal{F}$, if for all $\epsilon > 0$,
$$\mathbb{P}(\Tilde{F}_{\beta}(h^{(s)}, \alpha) - {F}_{\beta}(h_{\mathcal{F}}, \alpha) > \epsilon) \rightarrow 0 \text{ as } s \rightarrow \infty $$
where $h_{\mathcal{F}}$ is the best estimator function in the function space $\mathcal{F}$. 
\label{def1}
\end{definition}

For a fixed $h$, it is straightforward to see that the empirical feature-based risk is an unbiased and consistent estimate of the true feature-based risk. The consistency is proven by the law of large numbers in \cite{VS11}. To prove that the estimate is unbiased, we use the following equality 
\begin{eqnarray*} 
\begin{aligned}
\E(\tilde{F}_{\beta}(h^{(s)}, \alpha)) &= \E\bigg(\alpha + \sum_{t =1}^s \frac{[L(h(\mathbf{z}_t),d_t)-\alpha]^+}{ (1-\beta)s }\bigg)\\
&\rightarrow  \alpha + \frac{s}{(1-\beta)s} \E{[L(h,\tilde{d})-\alpha]^+}={F}_{\beta}(h_{\mathcal{F}}, \alpha).
\end{aligned}
\end{eqnarray*} 
as $s \rightarrow \infty$. The first and last equalities follow by the definitions of $\tilde{F}_{\beta}$ and ${F}_{\beta}$ respectively. The second equality follows from the linearity of expectation.  

\subsubsection{Function Space for $h$
\label{sub:fun_space}}

For a non-fixed function $h \in \mathcal{F}$, the function space selection is crucial to ensure consistency. Generally, the space should be large enough to capture the true estimator of $h$, but not so large that it becomes overly flexible and over-fits the randomness of the data.

To prove the uniform consistency, we use the Uniform Law of Large Numbers (ULLN). For a non-fixed function $h$,  the ULLN and Glivenko-Cantelli Theorem state that as the sample size increases, the empirical distribution converges uniformly to the true distribution (see \cite{VW96}, \cite{BM02}). This approach is often used in non-parametric estimation problems, where the true distribution is unknown and needs to be estimated from the data. To ensure consistency, we need the following assumptions. 

\begin{assumption} \label{assumption 2}
Assume the estimator $h$ and the function space $\mathcal{F}$ satisfies:
\begin{enumerate}
    \item The true estimator $h_{\mathcal{F}}(\cdot)$ exists and generates the data with random noise $d_t = h_{\mathcal{F}}(\boldsymbol{z}_t) + \Delta_t$, where $\Delta_t$ is a random variable with zero mean and finite variance.
     \item  The function space $\mathcal{F}$ is a Glivenko-Cantelli (`GC') class. (See \cite{Du02}, \cite{VW96}). 
\end{enumerate}
\end{assumption}

\begin{remark}
\label{remark: valid function space}
Using Assumption \ref{assumption 2}, the empirical risk converges uniformly to the true risk over the function space, $\mathbb{P}(\sup_{h \in \mathcal{F}} |\Tilde{F}_{\beta}(h^{(s)}, \alpha) - {F}_{\beta}(h_{\mathcal{F}}, \alpha) |) \rightarrow 0$ as $s \rightarrow \infty$.  

Assumption \ref{assumption 2} provides theoretical requirements for the function space. In practice, \cite{VW96} and \cite{Va98b} have demonstrated that the function class must satisfy two main conditions, namely, boundedness and equicontinuity. The following are some examples of function spaces that ensure consistency.

\begin{itemize}
    \item Lipschitz continuous functions: A class of real-valued functions on a compact domain where there exists a Lipschitz constant $L_p$ such that for all $h \in \mathcal{F}$ and ${z}_1, {z}_2$ are sample points, we have $\|h(\boldsymbol{z}_1) - h(\boldsymbol{z}_2)\| \le L_p  \|\boldsymbol{z}_1 - \boldsymbol{z}_2\|$ where $\|\cdot\|$ represents the Euclidean norm.  Lipschitz functions have a bounded rate of change, ensuring that their empirical process converges uniformly to the true process. For instance, \eqref{function h} is a linear function in a compact domain and thus is Lipschitz continuous. Also, polynomials in a bounded domain and log-linear functions on a compact domain excluding the origin are Lipschitz continuous. Details and proofs can be found in  Section 1.12 in \cite{VW96} and \cite{CC87}.      
    \item Functions with a finite Vapnik-Chervonenkis (VC) dimension. Yet, we note that within these functions, no specific examples stand out as particularly suited for our NPC method, as these functions are commonly used in classification problems.  Details can be found in Section 2.6 in \cite{VW96} and \cite{DGL96}. We included it here for the sake of completeness.
    \end{itemize}
    
\end{remark}
Lastly, we want to highlight that the GC class function serves as a theoretical foundation for understanding the uniform convergence of empirical risk to true risk in $\mathcal{F}$ with increasing sample size. However, it is important to note that the theoretical results represent the worst-case complexity which may be overly pessimistic (\cite{VS11}). Therefore, the choice of function space ultimately depends on the specific problem and the properties of the true estimation function. 

\subsubsection{Approximation error and sample complexity}

To conclude this section, we briefly discuss the approximation error and sample complexity. So far, we have shown that $ \hat{F}_{\beta}(h, \alpha) \rightarrow {F}_{\beta}(h_{\mathcal{F}}, \alpha)$, as sample size grows, and that the function space $\mathcal{F}$ can be selected intuitively. However, the ultimate question remains: Is it possible to find the best estimator $h_*$ among all functions? To achieve this, we need both consistencies in $\mathcal{F}$ and the convergence of approximation risks to $0$ as sample sizes increase. Namely
\[
\Tilde{F}_{\beta}(h^{(s)}, \alpha) - {F}_{\beta}(h_*, \alpha) =  \bigg[ \underbrace{ \Tilde{F}_{\beta}(h^{(s)}, \alpha) - {F}_{\beta}(h_{\mathcal{F}}, \alpha) }_{\text{estimation error}} \bigg] +  \bigg[ \underbrace{\Tilde{F}_{\beta}(h_{\mathcal{F}}, \alpha) - {F}_{\beta}(h_*, \alpha)}_{\text{approximation error}} \bigg] \rightarrow 0
\]
 as $s \rightarrow \infty$.  
The first term (also known as `variance') deals with uncertainty from random sampling when we estimate the best function $h$ in $\mathcal{F}$ using a finite sample. If our function space guarantees consistency, we can bound the estimation error. The second term (also known as `bias'), unaffected by random factors, arises from the model's inability to perfectly represent the true input-output relationship. 


Nested function spaces, $\{\mathcal{F}_s\}_{s \ge 0}$ as proposed by \cite{CZ07}, can be used for this purpose. For instance, we can select the function space as polynomial spaces with the degree gradually increasing with $s$.  Under assumptions on $h_*$,
the approximation error converges to zero as the degree increases. How to find the best approximation function in each function class $\mathcal{F}_s$ and how fast the best approximation converges to $h_*$ are questions that fall under approximation theory. This subject has been extensively studied (see \cite{Tr19}, \cite{SV90}, \cite{Ch82}). Concerns about having sufficient data points for polynomial coefficients in smaller sample sets can be addressed by sample complexity analysis (see \cite{AC20, AH20, ZN23}). 



\section{Baseline Experiment} \label{se:art}

In order to assess the performance of the proposed NPC method and to understand its
strengths and weakness, we conduct experiments with hypothetical data in \textit{R 4.2.1} with an Apple M1 Pro (2021) machine. In Subsection \ref{sub:setup}, we discuss the setup for our baseline experiment. In Subsection \ref{sub:exp1} the simplest case is studied, in which the profit function is linear. The case in which the profit function is nonlinear is discussed in Subsection \ref{sub:exp2}. 

\subsection{Experimental setup} \label{sub:setup}

For our baseline experiments, we suppose initially that there are 4 features related to the demand, each containing 500 observations (the cases with other numbers of features will be discussed later). We generate each feature from a seasonal ARIMA process, and we generate the demand from a linear model:
\begin{equation} \label{eq:demand}
    d_t := b_0 + b_1 z_t^1 + b_2 z_t^2 + b_3 z_t^3 + b_4 z_t^4 + r_t,
\end{equation}
where $z_t^p$ is the realisation of feature $p$ at time $t$, and $r_t$ is a realised error generated by an additive (weighted) mixture of \texttt{rnorm()}, \texttt{rlaplace()} and \texttt{rt()} functions. The choice of $b_p$ for the features, $\boldsymbol \phi$ and $\boldsymbol \theta$ for the ARIMA process, and other parameters for the generation of error terms are all selected randomly. We choose a seasonal ARIMA model since it is one of the most popular statistical models
in the literature \citep[for example, see][]{SBBKN16}. The detailed parameter values for the baseline setup and the demand series can be seen in Table \ref{tab:base_par} and Figure \ref{fig:base_series} in \ref{app:base}. The results of the experiment with other values of parameters will be discussed in Section \ref{se:ex_exp}. The scripts of all experiments have been made available on Github \citep{Li22}.

All experiments are performed on a rolling-origin basis with 1-step ahead order forecast \citep{Ta00}, in which we select fix \emph{origin size} (holdout sample size) to be 50, 100, 150, 200, 250 and 300, and the \emph{iteration number} (number of shifted) to be 50, 100, 150 and 200. For each pair of \emph{origin size} and \emph{iteration number}, we use the following
two quantities as measurements:
\begin{itemize}
    \item $\beta$-Downside Loss ($\beta$-$DL$) = $\frac{1}{n} \sum_{t=1}^{n}L_t$: This measures the average value of the largest $(1-\beta)$ cases of losses, where $n= \lceil (1-\beta)*iteration \, number \rceil$, and $L_t$ is ranked in descending order. It is desirable for this value to be as small as possible.
    \item Service Level ($SL$) = $\frac{1}{N}\sum_{t=1}^{N}\mathbb {I}{(x_t \geq d_t)}$: This measures the proportion of cases in which the demand is successfully fulfilled, where $N=iteration \, number$, and $\mathbb {I}(\cdot)$ is the indicator function. In the ideal situation, $SL$ should be as close to the target service level as possible.
\end{itemize}

In the proposed method (``NPC"), we use \texttt{optim()} function from \texttt{stats} package for \emph{R} and the Limited-memory Broyden-Fletcher-Goldfarb-Shanno algorithm (L-BFGS) for the estimation of parameters of the model. The L-BFGS algorithm has been shown to perform well in similar nonlinear programming tasks in an NVP context \citep{LN89,LLS22}.

Since $L$ can be a function of any level of complexity, minimising (\ref{eq:empirical}) is, in general, a continuous nonlinear optimisation problem. General-purpose algorithms for nonlinear optimisation are not guaranteed to converge to the global minima due to the lack of everywhere-differentiability. Fortunately, the experiments in Section \ref{se:art} and Section \ref{se:ex_exp} indicate that this does not cause serious problems.

To get around the scale issues, we consider two benchmark methods, with which we can compute the relative $\beta$-$DL$ and relative $SL$ (note that the ``instability" issue and the ``negativity" issue of relative measurement do not incur in our experiment):
\begin{itemize}
    \item The benchmark method - Sample weighted average (``SA"): With  historical demand $[d_1,\dots,d_s]$, the order quantity $x_{s+1}$ is set to be a weighted average of empirical quantiles under Equation (\ref{eq:CVaR}).
    \item The benchmark method - Under correctly specified model (``UM"): A method uses \texttt{lm()} function from \texttt{stats} package for \textit{R} to forecast the next period demand considering all features and all observations, and determines the order quantity with Equation (\ref{eq:CVaR}). We make sure the distribution of error term is correctly assumed.
\end{itemize}

Besides the benchmarks, three competing methods are also considered:
\begin{itemize}
    \item A non-featured method (``NF") that applies the \texttt{auto.arima()} function from \texttt{forecast} package for \textit{R}  to the demand series itself in the forecasting phase, and determines the order quantity with Equation (\ref{eq:CVaR}).
    \item A superquantile regression method (``SQR") that uses  \texttt{rq()} to determine the order quantity \citep{RRM14}. 
    \item A method (``PLM") that uses \texttt{lm()} to forecast but with only observations from the `worst' $2\times (1-\beta)$ proportion of scenarios (the term `worst' is defined in Subsection \ref{sub:method}).
\end{itemize}

In the baseline experiment, we consider a linear profit function
\begin{equation} \label{eq:exlinear}
\pi(x,d):= 20\,\min\{x,d\} - 8\,x +3\,[x-d]^{+} + 7\,[d-x]^{+},
\end{equation}
and a nonlinear profit function
\begin{equation} \label{eq:exnonlinear}
\begin{aligned}
\pi(x,d):= 20\,\min\{x,d\} &- 8\,x -4\,[x-d]^{+}\\
&+ 5\,\E[\min \{[x-d]^{+},u\}] - 0.01\,\big([d-x]^{+}\big)^2,
\end{aligned}
\end{equation}
where $u\sim \mathcal{N}(30,5^2)$. These settings are consistent with the work of \cite{LLS22}. The optimal service levels to maximise the expected profit
for these two functions are 0.35 and 0.62.

\subsection{The results of baseline experiment}

Here, we present the results from our experiment, where the parameters described in Subsection \ref{sub:setup} are used. 

To get some sense of the experimental procedure and the interpretation of the relative measurement, we first present a na\"{i}ve example in \ref{app:hist}, where only NPC method is considered, together with benchmarks. 

The relative $\beta$-$DL$ and relative $SL$ in this example can be calculated as:
\begin{equation}
    \text{relative} \; \beta \text{-} DL =  \frac{DL_{SA}-DL_{NPC}}{DL_{SA}-DL_{UM}} = 93\%
\end{equation}
\begin{equation}
    \text{relative} \; SL =  1 - \left|\frac{SL_{NPC}-SL_{UM}}{SL_{SA}-SL_{UM}}\right| = 20\%
\end{equation}
These values can be interpreted as: by using NPC, the decision makers will suffer $93\%$ less loss than in the case of using SA in the worst $5\%$ of scenarios, and the service level they achieve will be $20\%$ closer to the target service level. Here we note that the service level achieved by CVaR minimisation may lie far away from the service level achieved by expectation maximisation due to their different natures. In practice, a decision maker will need to balance the benefit from reducing downside loss and the harm of decreasing service level (as they are inseparable in most cases). 

\subsubsection{Linear profit function}
\label{sub:exp1}

\begin{figure}[ht]
\centering
\caption{Relative 95\%-$DL$ under linear profit function}
\label{fig:linear95}
\begin{subfigure}[b]{0.48\textwidth}
\centering
\includegraphics[width=0.9\textwidth]{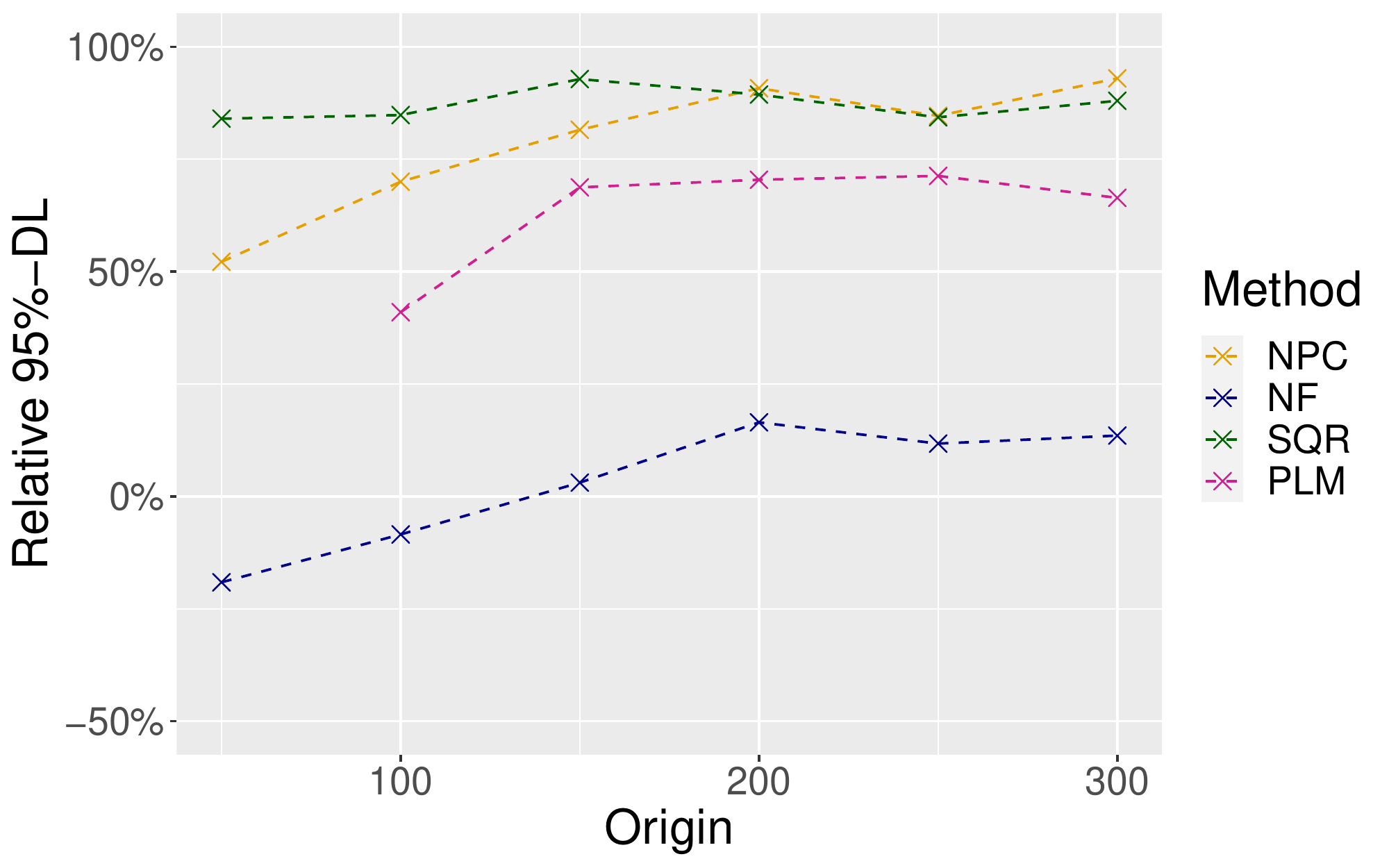}
\caption{Iteration number = 200}
\label{fig:draw_dl_x_200}
\end{subfigure}
\hfill
\begin{subfigure}[b]{0.48\textwidth}
\centering
\includegraphics[width=0.9\textwidth]{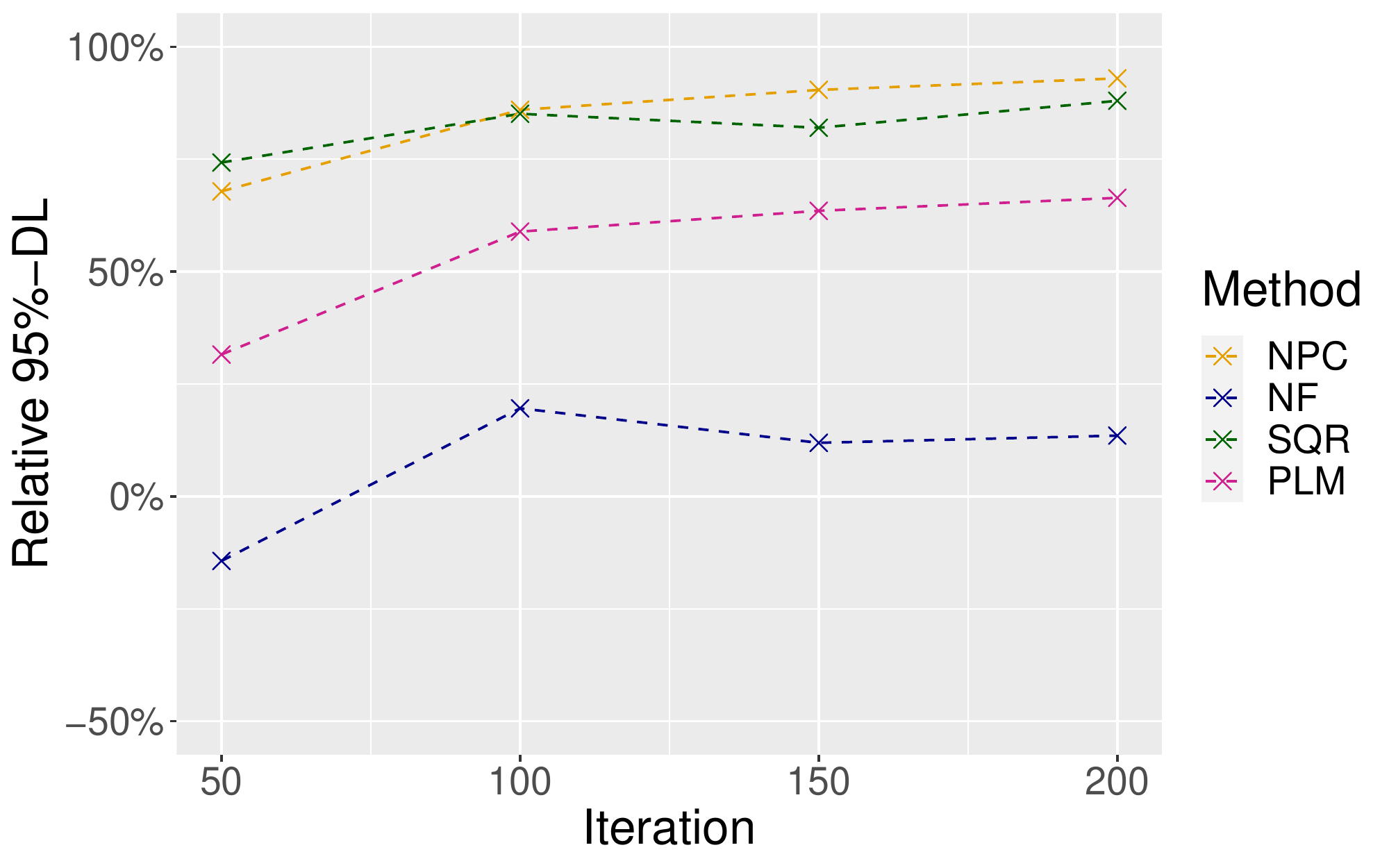}
\caption{Origin size = 300}
\label{fig:draw_dl_300_y}
\end{subfigure}
\end{figure}

In Figure \ref{fig:draw_dl_x_200}, we present the relative $\beta$-$DL$ from our baseline experiment, where multiple choices of \emph{origin size} are considered under the given linear profit function when \emph{iteration number} = 200. We remark that the result from PLM method under \emph{origin size} = 50 is excluded from the plot, as it is far below zero (This is probably due to the drawback of using Least-squares estimation under small sample size). In general, we could see that both the proposed NPC method and the SQR method achieve a high relative $\beta$-$DL$. In fact, their performance is quite close, even though SQR uses all historical observations and NPC only uses a small proportion of them. The results from the other two methods are less appealing, as PLM generates very frustrating performance when \emph{origin size} = 50, and NF barely improves the loss compared to the benchmark SA method.

In Figure \ref{fig:draw_dl_300_y}, we focus on the \emph{origin size} = 300 and present the relative $\beta$-$DL$ under multiple choices of \emph{iteration number}. We can see that the results are very similar to what we found in Figure \ref{fig:draw_dl_x_200}, where the NPC and SQR methods outperform the other two. In \ref{app:base_results}, we present the results from all other choices of \emph{origin size}, \emph{iteration number} and $\beta$ in detail, where we include the relative $SL$ as well.

We note that as we are using relative measurements, the results are seem to be ``stable" among all choices of parameters. This is to be expected, given that the absolute performance of all methods is influenced by parameters at the same time. From the results, we can say that the NFC method shows very strong robustness as its performance is very close to the SQR method (and the UM method) under all cases in regard to the relative $\beta$-$DL$. Using the same amount of data, the PLM method, however, performs poorly in most cases.

As there is no significant upward and/or downward trend in the demand series, as seen in Figure \ref{fig:base_series} in \ref{app:base}, it is totally understandable why the two non-featured methods, SA and NF, perform similarly. (Though the performance of NF improves slightly as the \emph{origin size} and \emph{iteration number} increase.) For the method of PLM, the nature of its loss function is to minimise the MSE, leading it to be under-fitted when the data is limited. With the same amount of data, the NPC method, on the other hand, adopts a different loss function and focuses on the extreme scenarios, making good use of all the selected data. The SQR method also performs well in this experiment. However, as we can see from Figures \ref{fig:draw_dl_x_200}, \ref{fig:draw_dl_300_y} and Table \ref{tab:linear_all} in \ref{app:base_results}, it gets slightly outperformed by NPC when \emph{origin size} and/or \emph{iteration number} is large, due to the presence of bias mentioned in Section \ref{se:intro}. This bias is amplified when the profit function is nonlinear and/or the error term distribution is changed, as we will see in Subsection \ref{sub:exp2} and \ref{sub:exp6}. 

\subsubsection{Nonlinear profit function} 
\label{sub:exp2}

Here, we present the results from our baseline experiment with nonlinear profit function. As one can see from (\ref{eq:exlinear}) and (\ref{eq:exnonlinear}), the major differences between these two forms of profit function are on the penalties of disposal and shortage. Instead of a fixed disposal cost, we now allow the excess products to be sold on a salvage market. Instead of a fixed shortage cost, we consider a quadratic cost function. As CVaR minimisation focuses on extreme cases, these differences may be amplified in our experiments and lead to results very different from those of Subsection \ref{sub:exp1}. Given that a closed form solution does not exist for the given nonlinear function, one can use the technique proposed by \cite{KK18}, or other numerical approaches, to verify that the quantiles to minimise 95\%-CVaR and 90\%-CVaR are approximately 0.13 and 0.16.

\begin{figure}[!ht]
\centering
\caption{Relative 95\%-$DL$ when iteration number = 200 under nonlinear profit function}
\label{fig:non_draw_dl_x_200}
\includegraphics[width=0.7\textwidth]{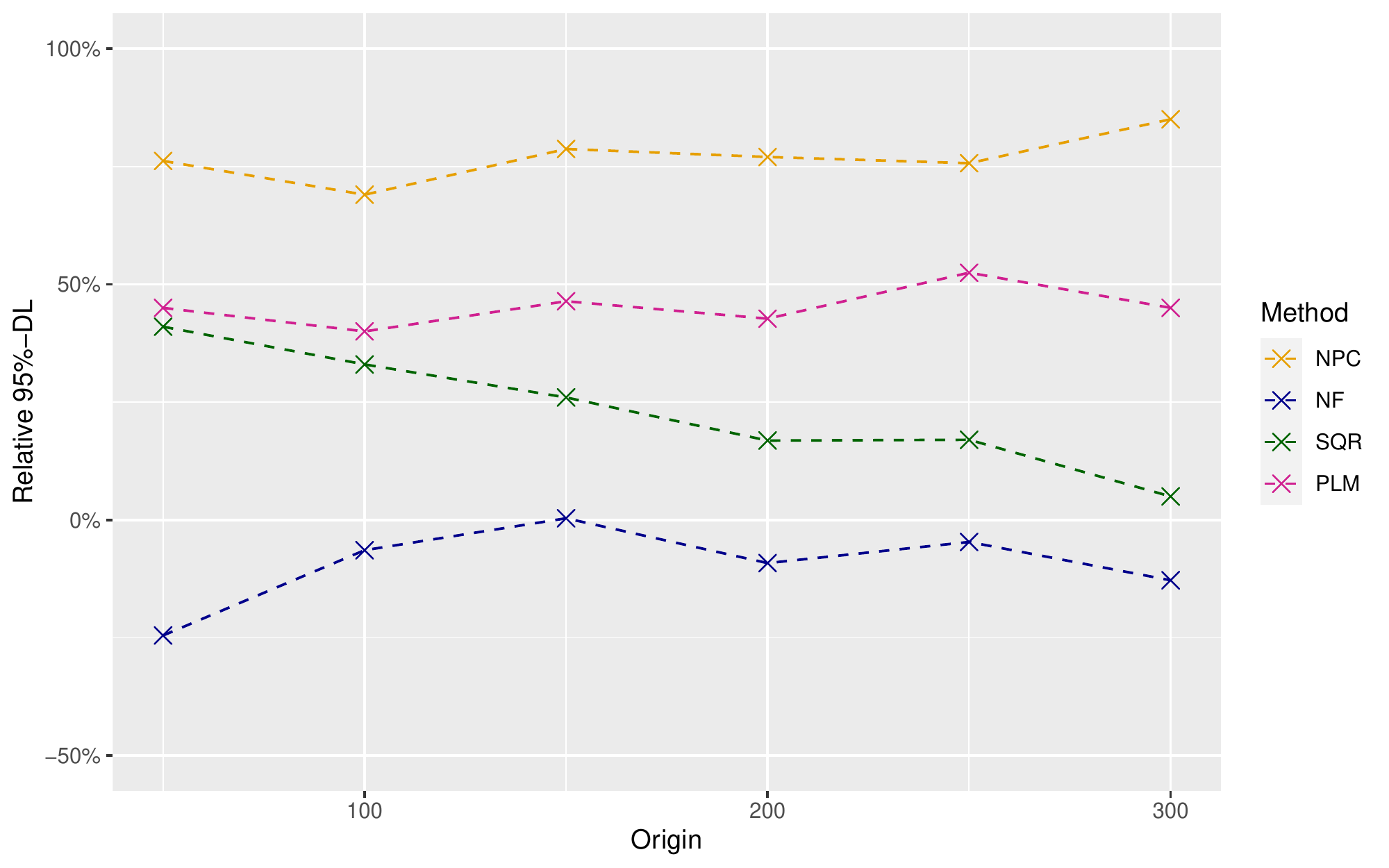}
\end{figure}

In Figure \ref{fig:non_draw_dl_x_200}, we present the relative $\beta$-$DL$ with multiple choices of \emph{origin size} when \emph{iteration number} = 200. It can be seen from the figure that the NPC method outperforms all other methods in regard to relative $\beta$-$DL$ under all \emph{origin size}. Moreover, we find that the relative performance of the SQR method decreases as \emph{origin size} increases. This can be further investigated by looking at the absolute performance in \ref{app:base_results}. We see that the $\beta$-$DL$ from SQR method does not improve as \emph{origin size} increases, while the $\beta$-$DL$ from SA does, leading to an overall decrease in relative $\beta$-$DL$. One possible explanation for this phenomenon is that in the SQR method, the loss function targets the extreme demand realisation instead of the extreme profit realisation directly. Therefore, under the nonlinear relationship between demand and profit, this loss function could be heavily biased. Thus, it is no surprise that the performance of SQR does not improve when increasing \emph{origin size}. On the other hand, the NPC targets the  extreme profit realisation.

We would like to stress that, unlike the parametric methods, NPC does not need any complicated numerical optimisation or simulation methods to estimate the optimal order quantity - it does that directly. In addition, NPC requires only a proportion of data under selection criterion, yielding results in a more efficient way. Overall, we see that NPC performs at least as well as SQR under linear profit functions, while outperforming all other methods under nonlinear profit functions. We will then examine the robustness of the NPC method in our next section.


\section{Experiments With Other Parameters}
\label{se:ex_exp}

Now, we extend our experiment to other parameters. In particular, we vary the numbers of features to be considered in Subsection \ref{sub:exp3}. 
Then, we present results with other profit functions in Subsection \ref{sub:exp5}. Final, we consider other forms of the error term in \ref{sub:exp6}. We remark that we have also experimented with other data generating models, e.g., ETS, TBATS. We do not present the results here as they are very similar to the ones presented below.

\subsection{Varying the number of features}
\label{sub:exp3}
Now, we focus on the number of features. In particular, we consider sufficiency of features adopted by the method (under-fitting/over-fitting), instead of the overall feature numbers, as this has negligible impact. This is motivated by the fact that in reality, decision makers are rarely able to guarantee the quality of feature choices \citep{HWD18}. Therefore, it makes sense for us to consider the performance of our proposed method in the case of model misspecification. To do that, we consider the relative $\beta$-$DL$ of the NPC method, the PLM method and one other method:
\begin{itemize}
    \item A regression method (``LM") that uses \texttt{lm()} to forecast with same number of features as used in NPC.
\end{itemize}
Besides, we also make sure the PLM method uses the same number of features as used in NPC and LM. We remark that the NF method and the SQR method are excluded from this comparison, for the obvious reason that they do not require any features in the computation. Without changing other settings, we now consider the cases where the method uses 3 features or 5 features instead, while using the same data set as before. These represent the cases of model under-fitting and model over-fitting, respectively. To avoid redundancy, here we only present the results with a linear profit function, as the results with a nonlinear profit function were very similar.

\begin{figure}[ht]
\centering
\caption{Relative 95\%-$DL$ when iteration number = 200 under linear profit function with other number of features}
\label{fig:features005}
\begin{subfigure}[b]{0.48\textwidth}
\centering
\includegraphics[width=0.9\textwidth]{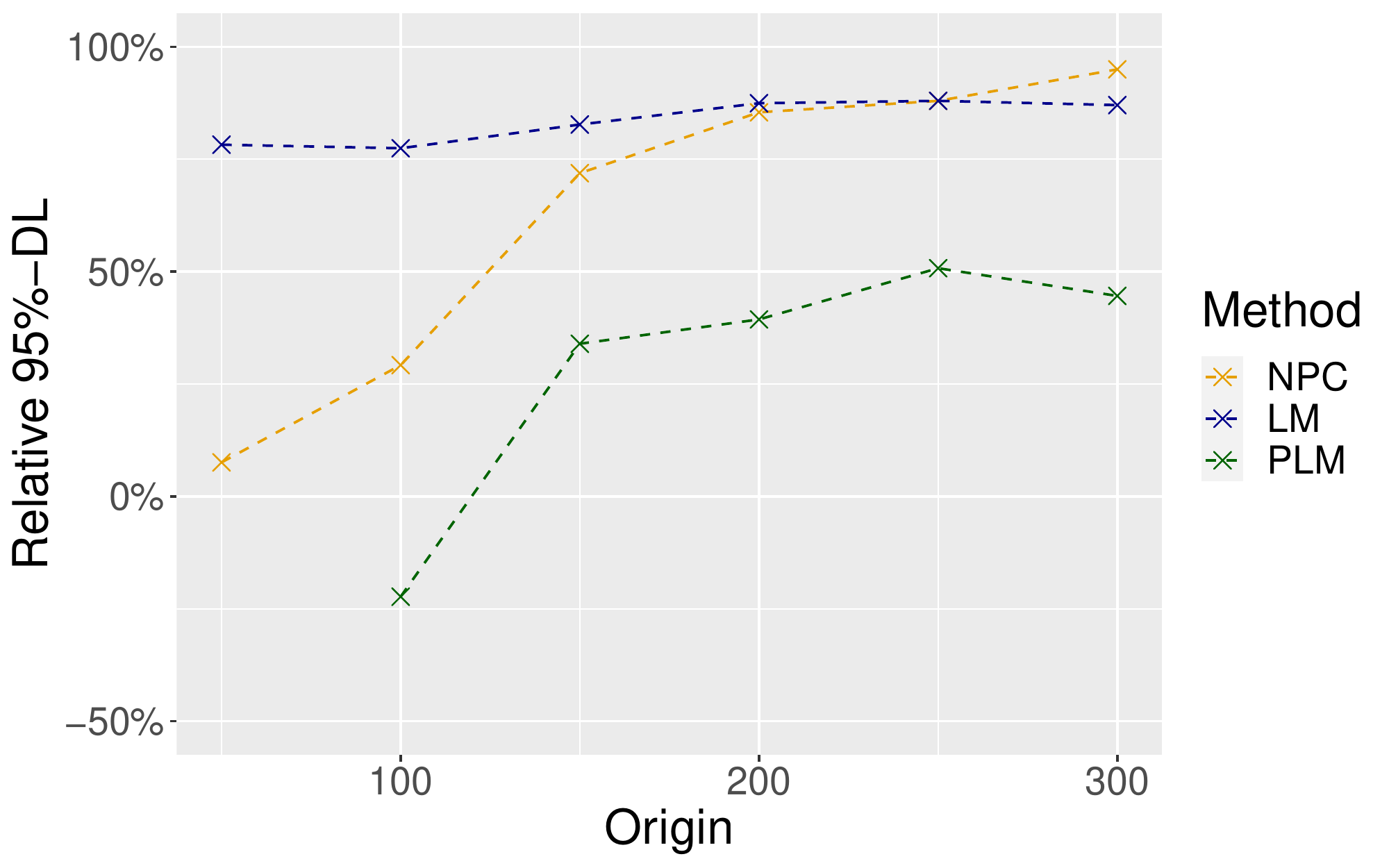}
\caption{With 3 features}
\end{subfigure}
\hfill
\begin{subfigure}[b]{0.48\textwidth}
\centering
\includegraphics[width=0.9\textwidth]{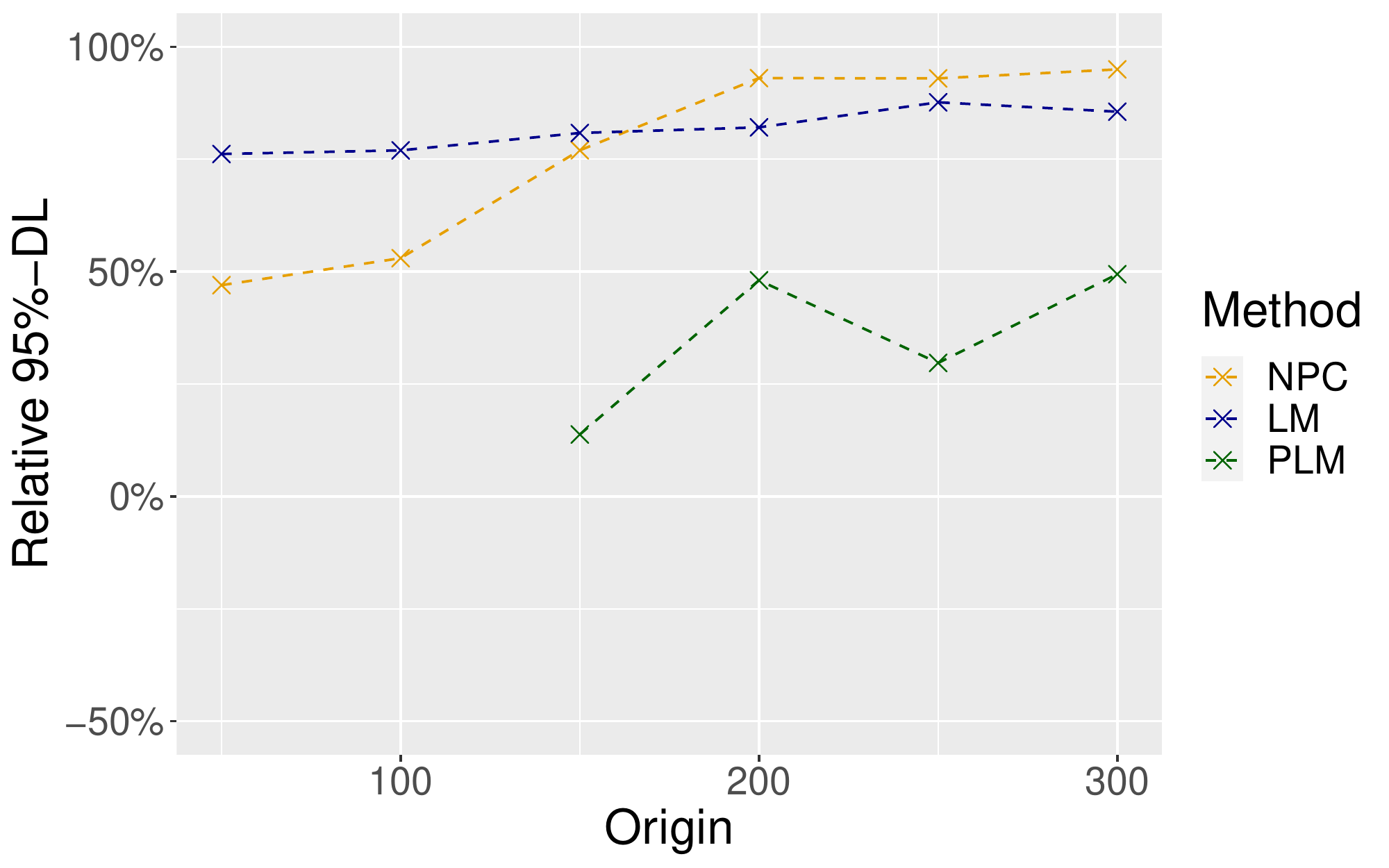}
\caption{With 5 features}
\end{subfigure}
\end{figure}

We can see from Figure \ref{fig:features005} that NPC performs better than PLM in all origin sizes, in both the under-fitting and over-fitting cases. However, its performance is worse than LM when \emph{origin size} is small, especially in the under-fitting case. The performance improves as \emph{origin size} increases. This is not completely unexpected. As the NPC method uses only a small proportion of the data, it could be more vulnerable than other methods when \emph{origin size} is small, especially when some information is missing due to under-fitting. Fortunately, we can see that, using the same amount of data, the performance of NPC is significantly better than the performance of PLM.

\begin{figure}[ht]
\centering
\caption{Relative 90\%-$DL$ when iteration number = 200 under linear profit function with other number of features}
\label{fig:features010}
\begin{subfigure}[b]{0.48\textwidth}
\centering
\includegraphics[width=0.9\textwidth]{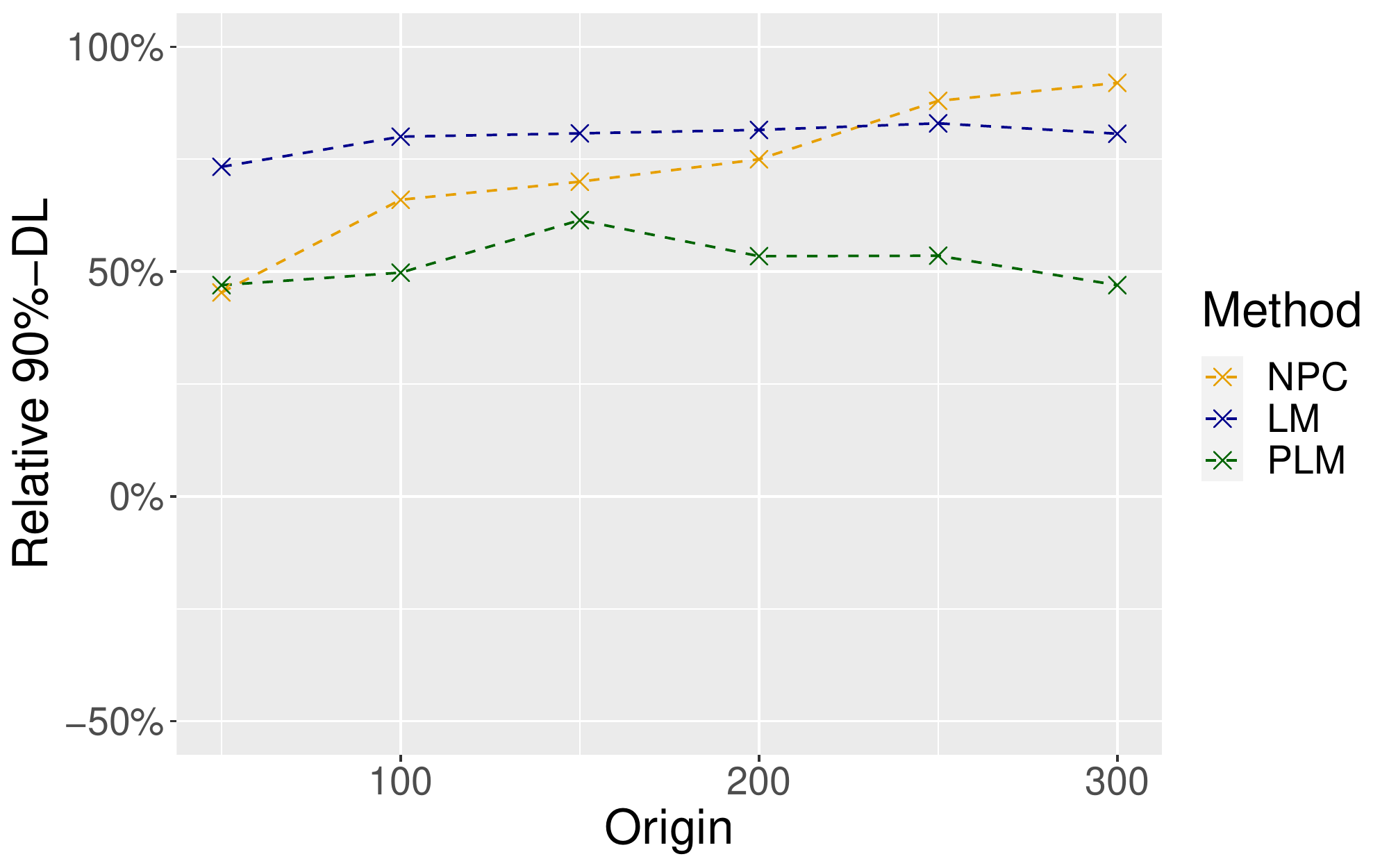}
\caption{With 3 features}
\end{subfigure}
\hfill
\begin{subfigure}[b]{0.48\textwidth}
\centering
\includegraphics[width=0.9\textwidth]{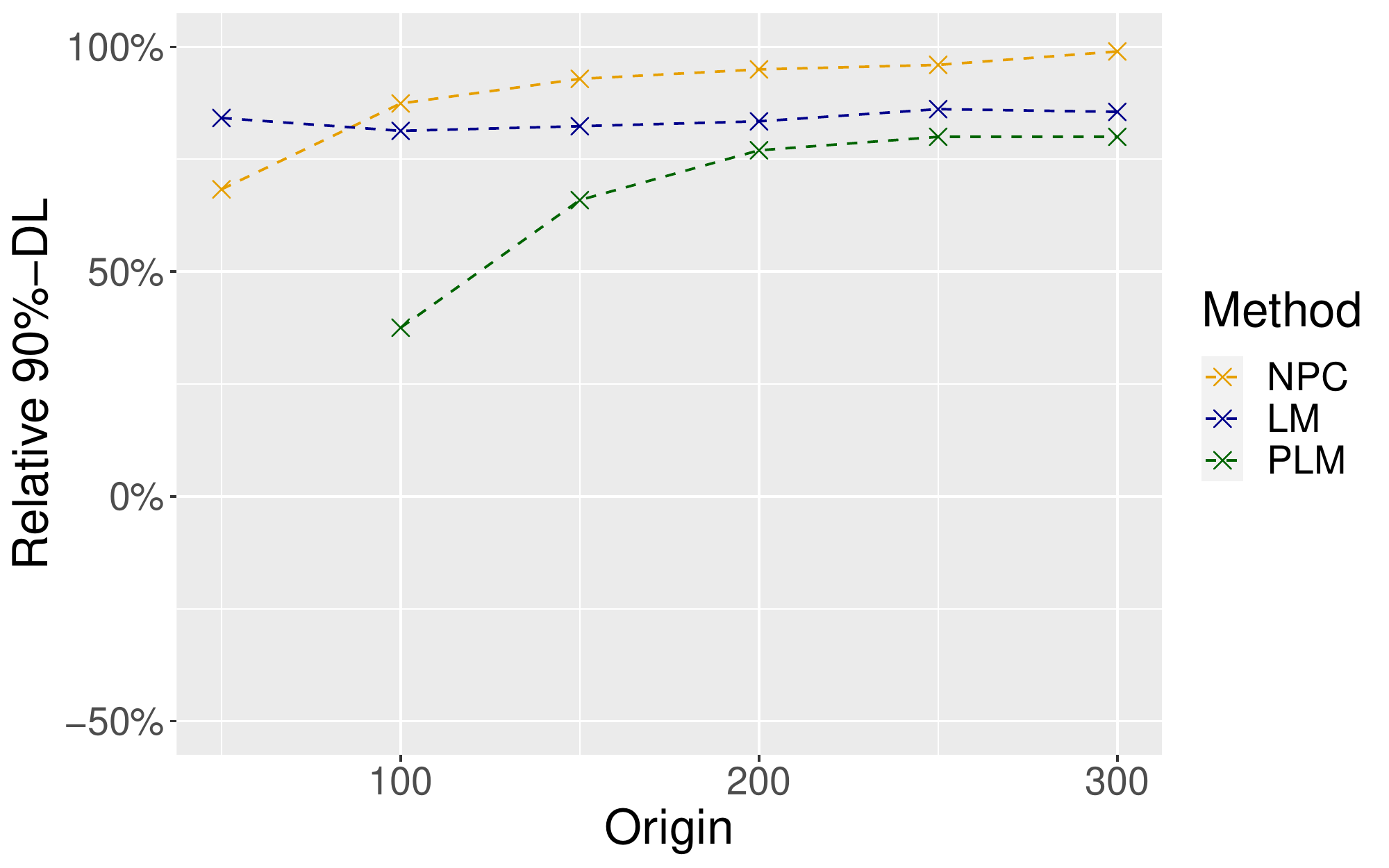}
\caption{With 5 features}
\end{subfigure}
\end{figure}

In Figure \ref{fig:features010}, we present the results where $\beta = 90\%$. In this setting, more data is used in the NPC method and PLM method. We can see that the performance of NPC is still slightly worse than LM when \emph{origin size} is small, but the gap is much smaller than in the case when $\beta = 95\%$. Besides, we find that the NPC method outperforms LM as long as the \emph{origin size} is larger than 100 in the over-fitting case, and 250 in the under-fitting case. We remark that the results with other \emph{iteration number} were very similar to the case when \emph{iteration number} = 200. Therefore, we do not present them here.

To sum up, we find that the proposed NPC method is more vulnerable than other methods when origin size is small, especially when the model is under-fitting. Nonetheless, this drawback is not unbearable, as our motivation in proposing an alternative method was to reduce the computational effort with large instances. Even in the case when \emph{origin size} = 300 (where NPC outperforms LM), the NPC method requires only 30 observations with $\beta = 95\%$, fewer than that required by LM when the \emph{origin size} = 50. 

\subsection{With other profit functions} 
\label{sub:exp5}

In the previous subsections, we tested the performance of our approach with one linear profit function and one nonlinear profit function, under different conditions. In this subsection, we consider four additional profit functions, two linear and two nonlinear, to examine the sensitivity of our method to the parameters of the profit function. All other settings are consistent with our baseline experiment.  We call the functions in the baseline experiment ``Linear 0" and ``Nonlinear 0", and we define ``Linear 1", ``Linear 2" and ``Nonlinear 1" and ``Nonlinear 2" as follows:
\begin{itemize}
    \item Linear 1:
    \begin{equation}
    \pi(x,d)= 20\,\min\{x,d\} - 8\,x -3\,[x-d]^{+} -7\,[d-x]^{+}.
    \label{eq:linear_1}
    \end{equation}    
    \item Linear 2:
    \begin{equation}
    \pi(x,d)= 20\,\min\{x,d\} - 8\,x +7\,[x-d]^{+} +3\,[d-x]^{+}.
    \end{equation}
    \item Nonlinear 1:
    \begin{equation}
    \pi(x,d)= 20\,\min\{x,d\} - 8\,x -4\,[x-d]^{+}- 0.01\,\big([d-x]^{+}\big)^2.
    \label{eq:nonlinear_1}
    \end{equation}
    \item Nonlinear 2:
    \begin{equation}
    \pi(x,d)= 20\,\min\{x,d\} - 8\,x + 5\,\E[\min \{[x-d]^{+},u\}],
    \end{equation}
    where $u\sim \mathcal{U}(0,15)$. 
\end{itemize}
``Linear 1" and ``Linear 2" are consistent with the work of \cite{LLS22}, while ``Nonlinear 1" and ``Nonlinear 2" are derived from it. The optimal service levels to maximise the expected profit for these four functions are 0.63, 0.9, 0.56 and 0.71, respectively. Although the order quantity that minimises the CVaR is usually very different from the one that achieves maximum expected profit, we find the optimal service levels of functions influence the performance of NPC.

\begin{table}[htb]
\caption{Relative 95\%-$DL$ under other profit functions (negative values are excluded)} \label{tab:pro_results}
\centering
\resizebox{0.8\linewidth}{!}{
\begin{tabular}{ccccc}
\hline
\multicolumn{1}{c}{} & \multicolumn{4}{c}{Methods}\\
\hline
Origin size = 50/Iteration number = 50 & NPC & NF & SQR & PLM\\
\hline
Linear 0 (0.35) & 60\% & / & $\mathbf{68\%}$ & /\\
Linear 1 (0.63) & 50\% & 3\% & $\mathbf{60\%}$ & 13\%\\
Linear 2 (0.9) & $\mathbf{66\%}$ & 28\% & $57\%$ & /\\
Nonlinear 0 (0.62) & $\mathbf{85\%}$ & / & 1\% & /\\
Nonlinear 1 (0.56) & $\mathbf{56\%}$ & 4\% & 5\% & 16\%\\
Nonlinear 2 (0.71) & $\mathbf{77\%}$ & / & 20\% & /\\
\hline
\multicolumn{1}{c}{} & \multicolumn{4}{c}{Methods}\\
\hline
Origin size = 300/Iteration number = 200 & NPC & NF & SQR & PLM\\
\hline
Linear 0 (0.35) & $\mathbf{90\%}$ & / & 72\% & 47\%\\
Linear 1 (0.63) & $\mathbf{84\%}$ & 13\% & 67\% & 30\%\\
Linear 2 (0.9) & $\mathbf{94\%}$ & 1\% & 74\% & 27\%\\
Nonlinear 0 (0.62) & $\mathbf{93\%}$ & / & 4\% & /\\
Nonlinear 1 (0.56) & $\mathbf{92\%}$ & 10\% & 49\% & 75\%\\
Nonlinear 2 (0.71) & $\mathbf{97\%}$ & 54\% & 39\% & 43\%\\
\hline
\end{tabular}}
\end{table}

As before, we present the Relative 95\%-$DL$ of all methods when \emph{origin size} = 50 and \emph{iteration number} = 50 with a linear profit function, as well as the case when \emph{origin size} = 300 and \emph{iteration number} = 200. We can see from Table \ref{tab:pro_results} that the NPC method perform well under all linear and nonlinear profit functions, and its performance converges to the UM method when \emph{origin size} and \emph{iteration number} increase. Specifically, the NPC method achieves better results when the optimal service level of the profit function is away from 0.5. This normally means that either the overage cost or the underage cost is much higher than the other. It appears that, when the profit function is heavy skewed, the NPC method is more efficient than competing methods to prevent downside loss in extreme cases. The phenomenon can be explained by the difference between the nature of the NPC method and competing methods. All competing methods compute results indirectly, as they work with extreme demand observations first and apply the output to the risk function second. However, the method of NPC works with extreme risks directly.

\subsection{With other forms of the error term} 
\label{sub:exp6}

Finally, we consider the influence of the error term. In our baseline experiment, the error term was generated by a mixture of \texttt{rnorm()}, \texttt{rlaplace()} and \texttt{rt()} functions with random parameters. Therefore, we have not yet examined how the proposed NPC method performs in
the presence of heavy-tails or light-tails. Given the drawbacks of traditional parametric methods on treating outliers, we could expect the gap of performance between NPC and PLM to be larger with light-tailed error terms than with heavy-tailed ones. We don't focus on the comparison of NPC and SQR in this case as they are influenced by the distribution of error term in a similar way. The experiment is conducted with some additional instances. We call the error term in the baseline experiment “Error 0”, and we define “Error 1” and “Error 2” as follows:
\begin{itemize}
    \item Error 1: We use \texttt{rnorm()} with $\mu = 0$ and $\sigma = 100$ as a light-tail case.
    \item Error 2: We use \texttt{rt()} with $\mu = 0$, $\sigma = 100$, $\nu = 5$ as a heavy-tail case.
\end{itemize}
We remark that it is not possible for the decision maker to know the exact distribution of the error term \emph{a priori} in reality. Therefore, in our experiment, we let our parametric methods assume that the distribution is normal in all cases, and we make sure that in our setting, the variance is the same in each instance.

\begin{table}[htb]
\caption{Relative 95\%-$DL$ under other forms of the error term (negative values are excluded)} \label{tab:error_results}
\centering
\resizebox{0.8\linewidth}{!}{
\begin{tabular}{ccccc}
\hline
\multicolumn{1}{c}{} & \multicolumn{4}{c}{Methods}\\
\hline
Origin size = 50/Iteration number = 50 & NPC & NF & SQR & PLM\\
\hline
Error 0 & 60\% & / & $\mathbf{68\%}$ & /\\
Error 1 & 71\% & / & $\mathbf{87\%}$ & /\\
Error 2 & 55\% & / & $\mathbf{82\%}$ & /\\
\hline
\multicolumn{1}{c}{} & \multicolumn{4}{c}{Methods}\\
\hline
Origin size = 300/Iteration number = 200 & NPC & NF & SQR & PLM\\
\hline
Error 0 & $\mathbf{86\%}$ & / & 72\% & 49\%\\
Error 1 & $\mathbf{90\%}$ & 6\% & $85\%$ & 40\%\\
Error 2 & $\mathbf{84\%}$ & 9\% & 84\% & 60\%\\
\hline
\end{tabular}}
\end{table}
The results in Table \ref{tab:error_results} meet our expectation, as the gap of performance between NPC and PLM is indeed very large in ``Error 1", and it is smaller in ``Error 2". A possible explanation is that the NPC method works directly with the data, and does not rely on the assumption of normality, while
the PLM, using the same amount of data, relies on normality. The gap is largest under light tails, since that the downside loss in extreme cases in light tail distribution is more likely to be treated as outliers than in heavy tail by traditional methods like PLM. Moreover, considering NPC itself, its performance improves as data size becomes larger and it performs best against other methods when the distribution of the error term is light tailed. Although SQR performs slightly better than NPC when data size is small, the NPC outperforms again once data size becomes large.

\subsection{Sum up}

From the results of experiments, we find that:
\begin{enumerate}
    \item Comparing to competing methods, NPC performs better when the data size is large. Besides, due to the fact NPC only require a small proportion of data, it can significantly reduce computational effort under large data set.
    \item The NPC performs well with nonlinear NVP, since it does not make any assumption on the linearity of the functions and works directly on historical data.
    \item NPC is more vulnerable than other methods when data size is small, especially when some information is missing due to under-fitting. 
    \item The proposed method performs best, comparing to competing method, when the profit function is heavy skewed and/or when the distribution of the error term is light tailed. In both cases, the traditional methods are likely to treat extreme data as outliers and underestimate the downside risk.
\end{enumerate}


\section{Real-life Example} \label{se:exp_real}

In this section, we examine the performance of NPC with a real-life example within a food bank. 
A food bank is an emergency feeding organisation providing hunger relief to families living in poverty. Each food bank covers a given region, and the decision maker has to prepare food on weekly basis for its distribution day (normally on Sunday). This food preparation problem within food bank can be approximately fitted by the nonlinear NVP model. The goal of the problem is to determine the amount of food to prepare that fulfils the demand. In the simplest case, we assume the consumption of each individual is same, and we could just use the 'number of visit' as our demand. Moreover, we assume both $x$ and $\tilde{d}$ (under same scale) to be continuous. Yet, we should note in particular that:
\begin{enumerate}
    \item The demand in food bank normally has smaller variance than the demand considered in other classic inventory management problem. Thus, instead of the expected profit, the CVaR is more of our interest.
    \item The opportunity cost of overage is linear since the food bank can easily dispose of excess food. In contrast, the cost of underage is thought to be quadratic, as the negative effects of not having enough food become increasingly severe and harder to mitigate over time. Moreover, many of these costs are intangible and difficult to quantify, which exacerbates their impact.
\end{enumerate}
Derived from \cite{DSIBM14} and \cite{Ri18}, this problem can be approximated as:
\begin{equation}
    L(x,d)= \eta[x-d]^{+} + \zeta \, \big( [d-x]^{+} \big)^2.
\end{equation}
where $\eta$ denotes the overage opportunity cost, including but not limited to transportation fee, management cost and disposal fee; $\zeta$ denotes the underage opportunity cost, including but not limited to loss of goodwill and additional management cost. The objective
is to minimise $L$, which we call `unit of risk' for generality. The parameter values are approximately chosen as:
\begin{itemize}
    \item $\eta = 15$, $\zeta = 1$
\end{itemize}
The parameter selection implies that the opportunity cost of underage increases rapidly as the supply-demand gap widens, eventually exceeding the opportunity cost of overage. In contrast, when the gap is small, the opportunity cost of overage is higher than that of underage. This is expected in food bank operations since small food shortages typically have only short-term effects, while large food shortages may damage the organisation's long-term reputation.

The data we use comes from a local food bank in Durham. It includes the total visit on each distribution day for 104 weeks from July 2020 to June 2022 on a weekly basis. Besides, we consider 10 relevant features within same time scope, as seen in Table \ref{tab:food_features} in \ref{app:invexp}.

\begin{figure}[!htb]
\centering
\caption{Time-series plot for the number of visit}
\label{fig:food}
\includegraphics[width=0.7\textwidth]{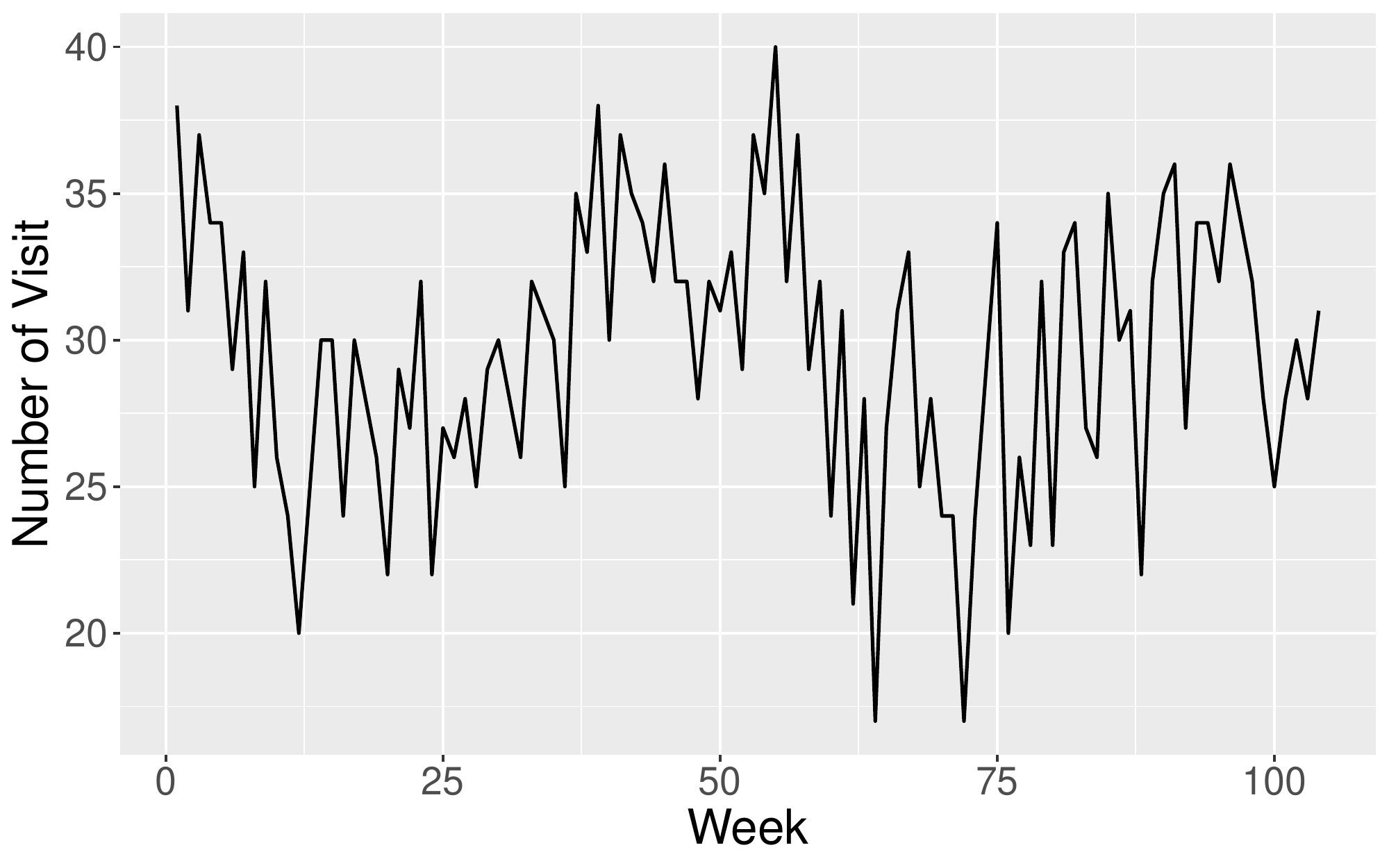}
\end{figure}

To get some sense of the data, we provide a time-series plot for the number of visit in Figure \ref{fig:food}.
It can be seen that the number of visit to food bank shows multiple levels of seasonality, monthly and seasonally, and that, rather surprising, the number of visit in winter (week 10-30 and 60-80) is lower than the rest of the year. We think this could probably due to the substitute effect from other forms of winter-exclusive aid, e.g., winter appeal, Christmas grants. 

Again, we use SA and UM methods as benchmarks. This time, we consider 10 methods that include different number of features: \begin{itemize}
    \item Non-feature: NF
    \item Seasonal feature (9-10): PLM-0, LM-0, NPC-0
    \item Local feature (5-10): PLM-1, LM-1, NPC-1
    \item National feature (1-10):
    PLM-2, LM-2, NPC-2
\end{itemize} 

\begin{table}[htb]
\caption{Relative performance when $\beta=0.95/0.90$ with 10 methods}
\label{tab:food_results}
\centering
\resizebox{0.8\linewidth}{!}{
\begin{tabular}{cccccc}
\hline
\multicolumn{1}{c}{$\beta$ = 0.95} & \multicolumn{5}{c}{Measurements}\\
\hline
Methods & rMAE & rMPS & rRMSE & Relative 95\%-$DL$ & Relative $SL$\\
\hline
NF & 91\% & / & 91\%  & 7\% & 0\%\\
PLM-0 & / & / & / & / & /\\
PLM-1 & / & / & 91\% & / & /\\
PLM-2 & / & / & 77\% & / & 10\%\\
LM-0 & 74\% & / & 78\%  & 2\% & 10\%\\
LM-1 & 39\% & 92\% & 65\% & 45\% & 15\%\\
LM-2 & $\mathbf{29\%}$ & 85\% & $\mathbf{50\%}$  & 92\% & 20\%\\
NPC-0 & 83\% & 88\%  & 72\% & 68\% & 15\%\\
NPC-1 & 49\% & 63\%  & 54\% & 85\% & 30\%\\
NPC-2 & $\mathbf{29\%}$ &$\mathbf{33\%}$ & $\mathbf{50\%}$ & $\mathbf{96\%}$ & $\mathbf{30\%}$\\
\hline
\multicolumn{1}{c}{$\beta$ = 0.90} & \multicolumn{5}{c}{Measurements}\\
\hline
Methods & rMAE & rMPS & rRMSE & Relative 90\%-$DL$ & Relative $SL$\\
\hline
NF & 89\% & 99\% & 97\%  & 13\% & 1\%\\
PLM-0 & / & / & / & 1\% & /\\
PLM-1 & / & / & 90\% & 1\% & /\\
PLM-2 & 92\% & / & 78\% & 5\% & 10\%\\
LM-0 & 70\% & / & 78\% & 4\% & 10\%\\
LM-1 & 38\% & 91\% & 66\% & 46\% & 15\%\\
LM-2 & 38\% & 85\% & 55\% & 93\% & 25\%\\
NPC-0 & 83\% & 89\%  & 72\% & 68\% & 15\%\\
NPC-1 & 48\% & 63\%  & 54\% & 88\% & 30\%\\
NPC-2 & $\mathbf{27\%}$ & $\mathbf{33\%}$ & $\mathbf{50\%}$ & $\mathbf{98\%}$ & $\mathbf{30\%}$\\
\hline
\end{tabular}}
\end{table}
To compare the performance of the methods, we obtain their 1-step ahead forecasts with rolling horizon, where \emph{origin size} is 60 and the origin is shifted 44 times. For each forecasted value, we compute the overage/underage amount and the cost. We summarise the results in Table \ref{tab:food_results}, where rMAE denotes the Relative Mean Absolute Error for the visit estimation, rMPS denotes the Relative Mean Pinball Score and rRMSE denotes the Relative Root Mean Square Error \citep{DF13}. We recall that low rMAE, rMPS and rRMSE are favourable, while high Relative $DL$ and Relative $SL$ are favourable.

From Table \ref{tab:food_results}, we can see that the NPC method with national feature performs best, having lowest error and highest relative DL. Moreover, under same number of features and same choice of $\beta$, the NPC outperforms other competing methods. This is expected. Since the true distributions of the time series and the error term are both unknown, and the cost function is nonlinear, the tail of the downside loss is hard to capture. All competing methods suffer from the underestimation of the downside loss from the tail.  Moreover, NPC requires only a proportion of input data, significantly improving the computing speed. For instance, when $\beta = 0.90$ and national feature is considered, the NPC method computes 5 times faster than the LM-2 method. Thus, NPC can not only help the decision maker to achieve lower downside risk, but also works more efficiently overall.

\begin{figure}[!htb]
\centering
\caption{Results of NPC-2 for week 61-104}
\label{fig:foodvs}
\begin{subfigure}[b]{0.48\textwidth}
\centering
\includegraphics[width=0.9\textwidth]{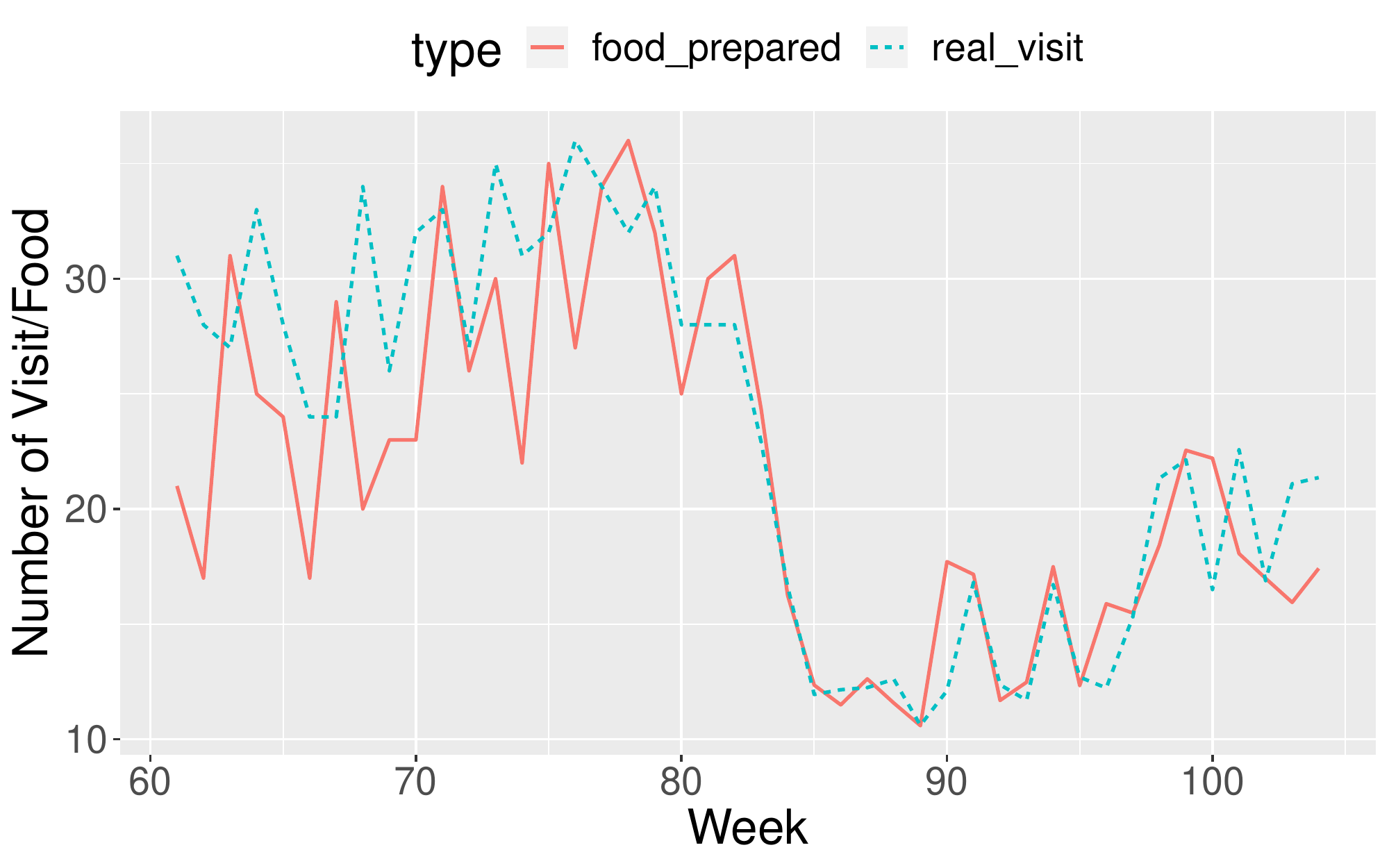}
\caption{Number of visit vs. food prepared}
\label{fig:fvv}
\end{subfigure}
\hfill
\begin{subfigure}[b]{0.48\textwidth}
\centering
\includegraphics[width=0.9\textwidth]{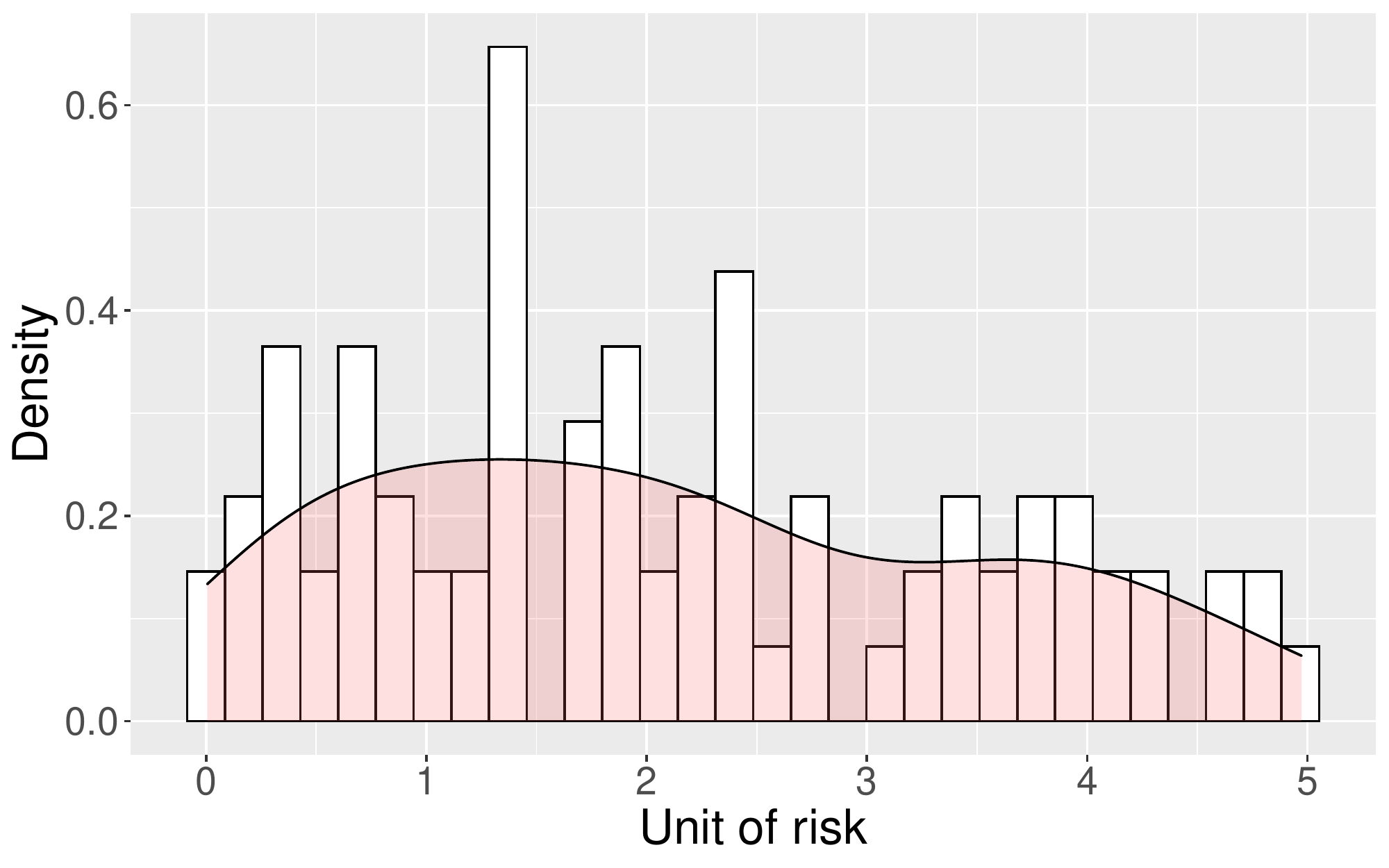}
\caption{Histogram of risk distribution}
\label{fig:his}
\end{subfigure}
\end{figure}

Last but not least, let us now suppose the food bank indeed implements our NPC-2 approach with $\beta =0.90$ for its food preparation  decisions. We wish to gain some insights into the predictions made by the approach. Figure \ref{fig:foodvs} provides the results acquired by NPC-2.

We derive the following insights from Figure \ref{fig:food} and {Figure} \ref{fig:foodvs}, which could be useful for both the visitor and the food bank:
\begin{enumerate}
    \item Due to the policy of minimising CVaR, the food bank may under-stock their food in winter. 
    \item For food bank visitors, if possible, we would suggest visit the food bank on summer, while try acquiring support from other source in winter.
    \item The risk generated using NPC-2 is close to `evenly' distributed according to Figure \ref{fig:his}. The food bank will not face significant loss, but it also cannot achieve minimum overall risk.
    \item The two most important features of food preparation are `unemployment rate' and `Covid-19 cases'. The food bank should increase its food preparation by 0.08 and 0.11 units as one unit of those two features increase, respectively.
    \item We suggest the food bank systematically obtain and record data that may be associated with the demand, since the proposed method performs better as more features are included. 
    \item Keep track of the performance of the model over time, as the parameters of model may evolve.
\end{enumerate}


\section{Concluding Remarks} \label{se:conclusion}

In this paper, we proposed an alternative non-parametric method (NPC) of CVaR minimisation with feature-based demand data. Unlike common non-parametric methods, the NPC method uses an adaptive data selection criterion and requires only a small proportion of data, significantly reducing the computational effort. Our proposed method works directly with the data, requiring no prior knowledge of the demand distribution. Additionally,  the estimated parameters of NPC can be easily applied to prescriptive analytic to provide additional operational insights. Our experiments with both hypothetical data and real-life data indicate that our proposed method is very robust with regards to different data structure, and it can easily handle both linear profits and nonlinear profits. In particular, one should be careful using NPC when the sample size is small, especially when the model is under-fitting, as it can be more vulnerable than other competing methods in this case. Yet, this drawback is not unbearable, as our motivation in proposing NPC was to reduce the computational effort with large instances.

There are several interesting topics for further research. First, as observed in our experiment, the performance of NPC suffers from the model under-fitting. Therefore, it would be interesting to extend the current NPC model to deal with this drawback. For instance, one can try introducing an additional parameter that controls the data usage manually (to a value other than $2 \times (1-\beta)$). Second, it would be desirable to develop a variable selection mechanism in NPC, as to prevent the model from over-fitting automatically, e.g., by cross-validation or a step-wise technique based on an information criteria. Finally, though we focused our research on NVP, the proposed method could be valuable in fields other than inventory control, such as in finance and logistics. Specifically, the convergence result we obtained can be applied to similar problems, provided that the objective function of the problem and the corresponding variables satisfy the assumptions made in Assumption \ref{assumption 1}. Such problems include, but are not limited to, portfolio optimisation problems under downside risk constraints, optimal product design problems under uncertainty, and supply chain risk management problems.
\bigbreak 
\bigbreak 
\noindent
\textbf{Acknowledgements}: The work was supported by the Hong Kong Innovation and Technology Commission (InnoHK Project CIMDA).


\newpage
\appendix

\section{Expectation maximisation vs. CVaR minimisation} \label{app:e_CVaR}

\begin{figure}[h]
\centering
\caption{Difference between the expectation maximisation solution and the CVaR minimisation solution}
\label{fig:mean_CVaR}
\includegraphics[width=0.7\textwidth]{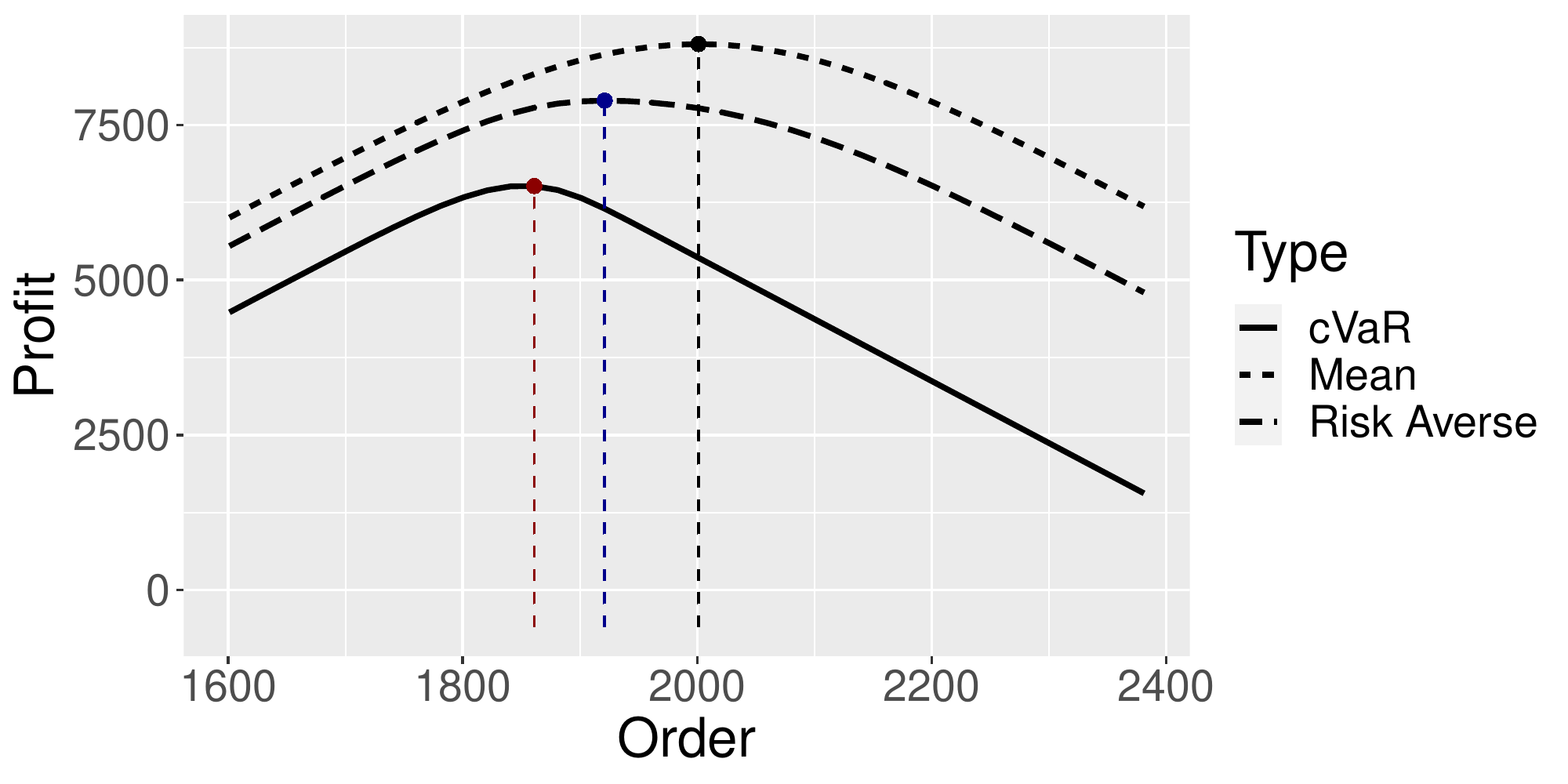}
\end{figure}

In Figure \ref{fig:mean_CVaR}, we mark three order quantities for a simple example where $E=10$ and $U=W=5$ with $\mathcal{N}(2000,150)$. They fulfil the objectives of expectation maximisation, CVaR minimisation and risk averse profit maximisation ($0.7\times \text{Mean} - 0.3\times \text{CVaR}$), respectively. We can see the order quantity that minimise CVaR is lower than the order quantity that maximise expectation. However, this is parameter-dependant, as the CVaR minimisation order quantity is a weight average of critical quantitles. When the overage cost is significantly larger than the underage cost, the CVaR minimisation quantity will be, with no doubt, larger than the expectation maximisation quantity. We argue that this result holds generally on other choices of parameters, referring to \eqref{eq:simCVaR}.

\newpage
\section{Numerical example for relationship between $S$ and $M$}
\label{app:sm}

\begin{figure}[!ht]
    \centering
    \caption{Illustrative example for Theorem \ref{thm2}. The linear loss function as in Equation \eqref{eq:linear_1} and a nonlinear loss function as in Equation \eqref{eq:nonlinear_1}. }
    \label{fig:res_distribution}
    \begin{subfigure}[b]{0.48\textwidth}
    \centering
    \includegraphics[width=0.8\textwidth]{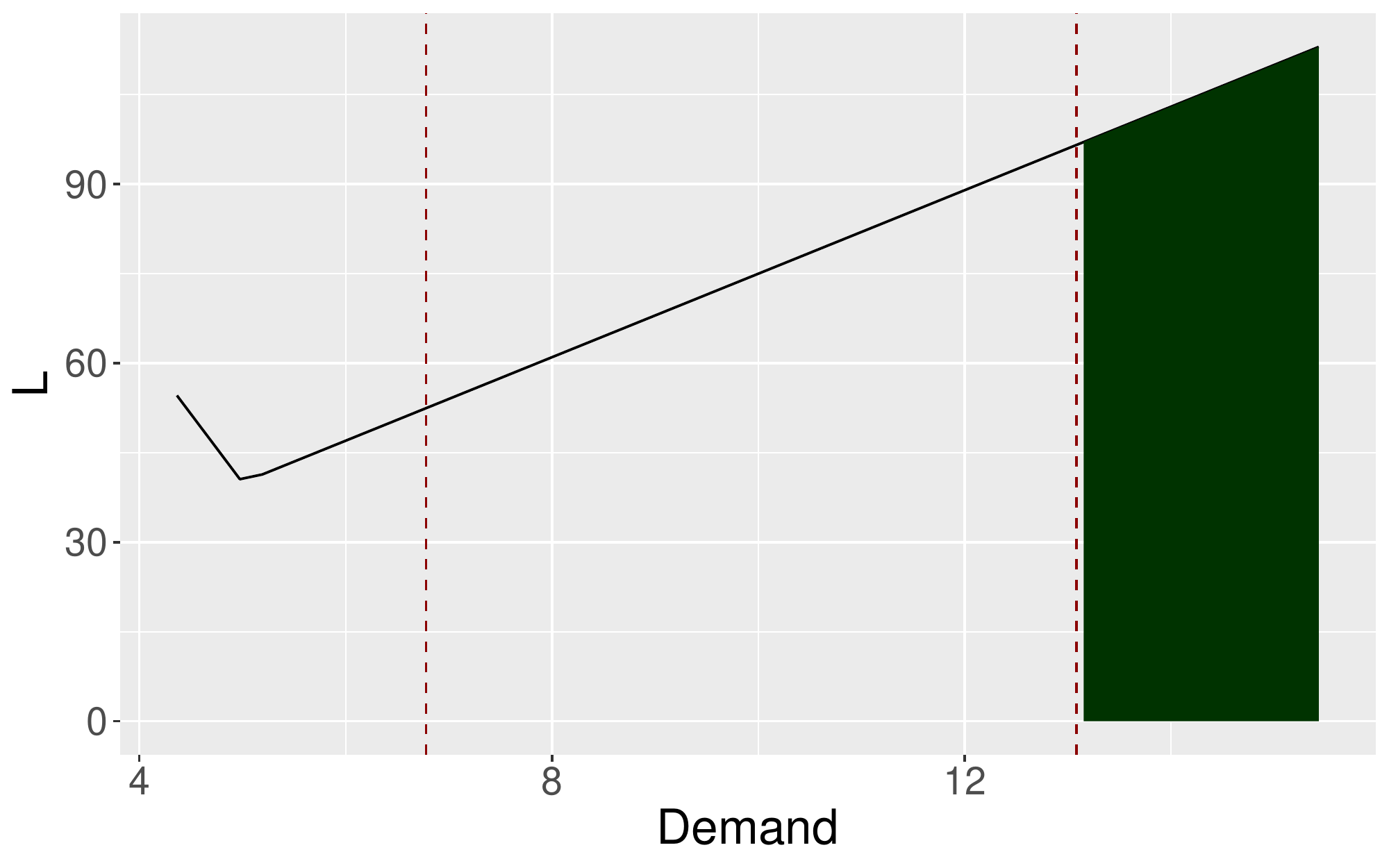}
    \caption{Linear loss function, $x=5$}
    \end{subfigure}
    \hfill
    \begin{subfigure}[b]{0.48\textwidth}
    \centering
    \includegraphics[width=0.8\textwidth]{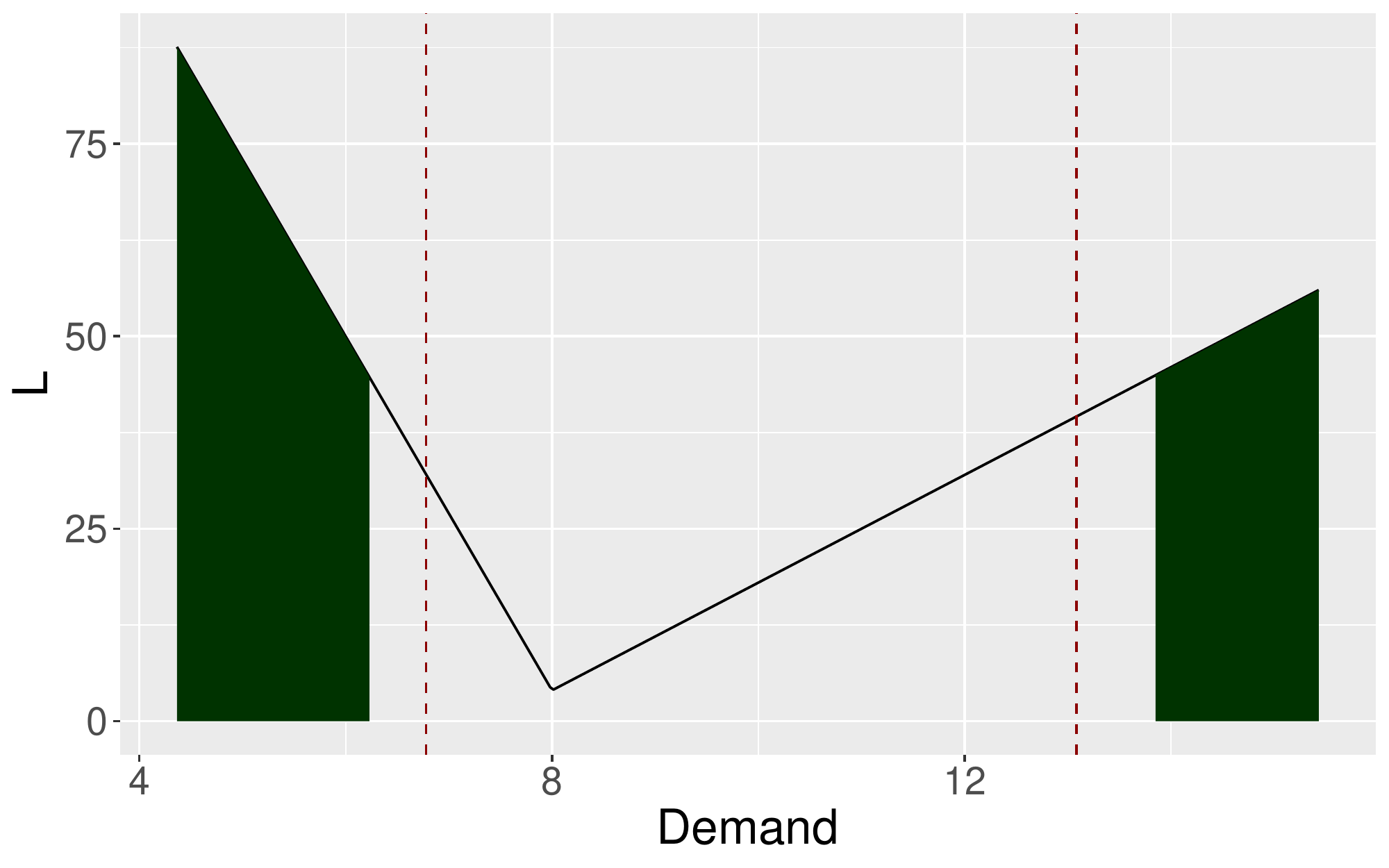}
    \caption{Linear loss function, $x=8$}
    \end{subfigure}
    \\
    \begin{subfigure}[b]{0.48\textwidth}
    \centering
    \includegraphics[width=0.8\textwidth]{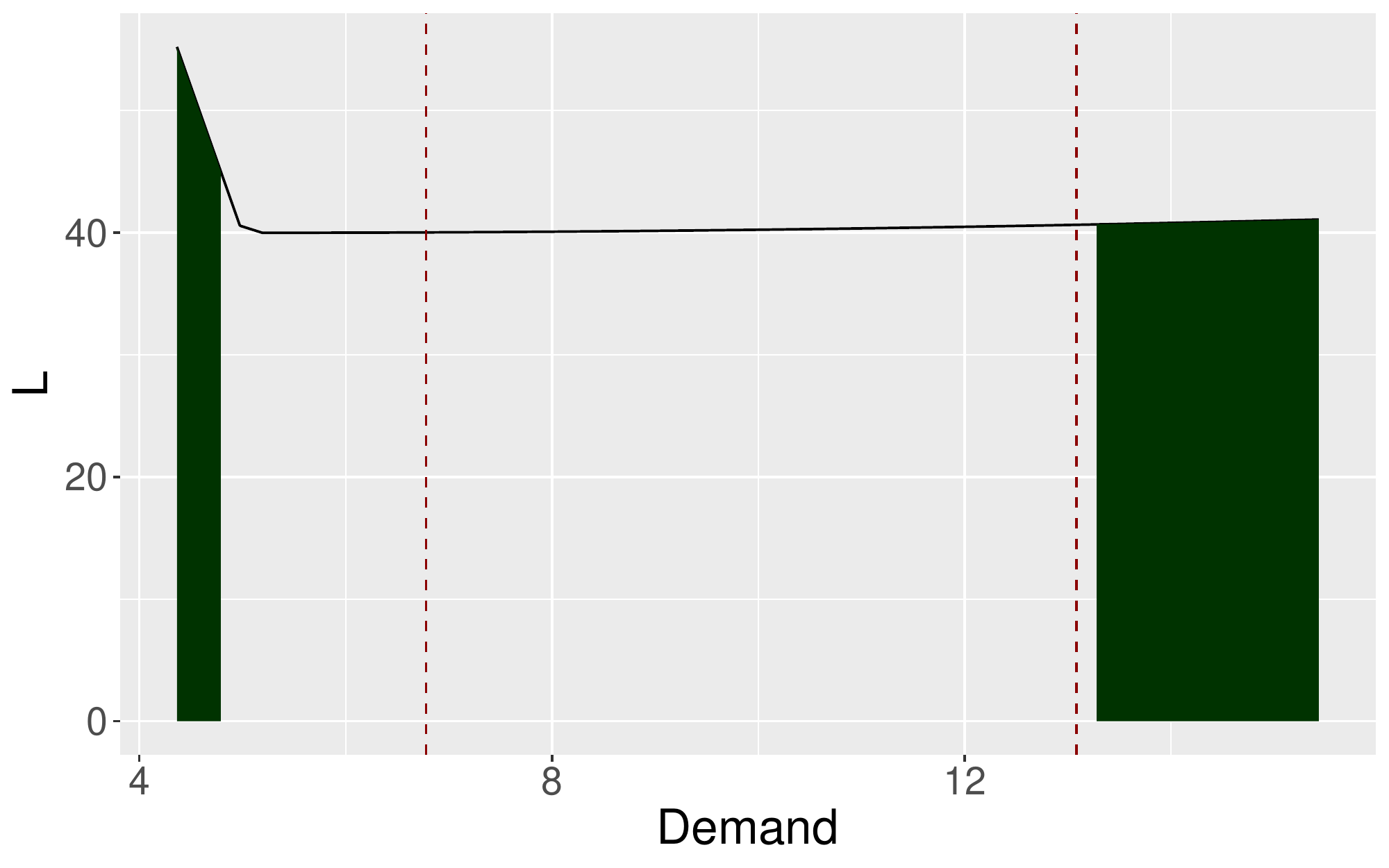}
    \caption{Nonlinear loss function, $x=5$}
    \end{subfigure}
    \hfill
    \begin{subfigure}[b]{0.48\textwidth}
    \centering
    \includegraphics[width=0.8\textwidth]{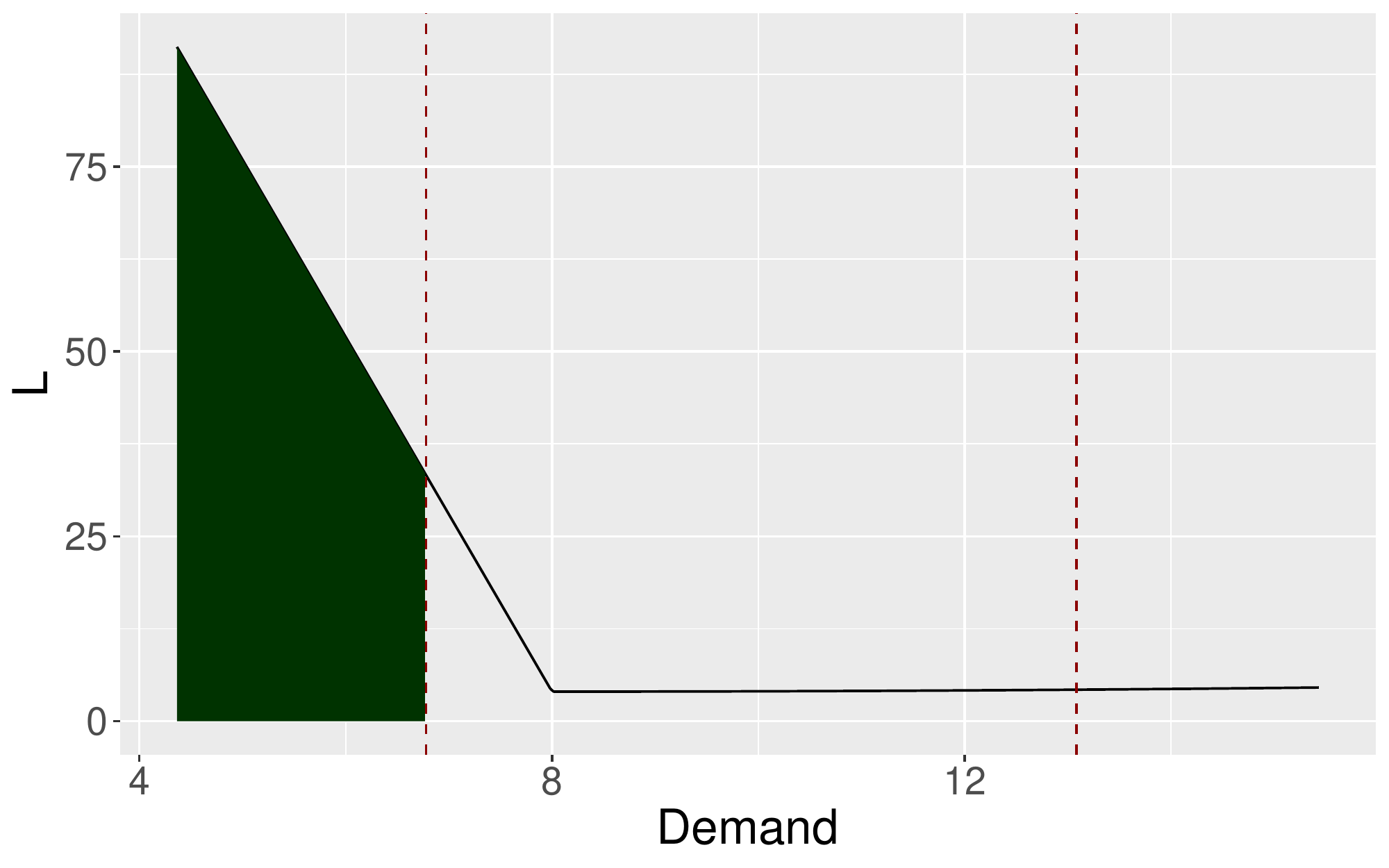}
    \caption{Nonlinear loss function, $x=8$}
    \end{subfigure}
\end{figure}

We illustrate the relationship between $S$ and $M$ with some numerical examples.  We generate the demands $d_t$ with a constant $T_t$ and an irregular component $\epsilon_t$ following a mean zero normal distribution. It is straightforward to verify that such distribution satisfies the distribution assumption. We consider a linear loss function and a nonlinear loss function in Figure \ref{fig:res_distribution}. Both loss functions satisfy the continuity assumption and the tail assumption. As a result, we see that for all plots, the set $S$ (shaded area) is contained in the set $M$ (the region bounded by dash lines and the vertical edges of the graphs).

\section{Baseline experiment parameters} \label{app:base}

\begin{figure}[htb]
\centering
\caption{Baseline experiment features and demand series}
\label{fig:base_series}
\includegraphics[width=0.7\textwidth]{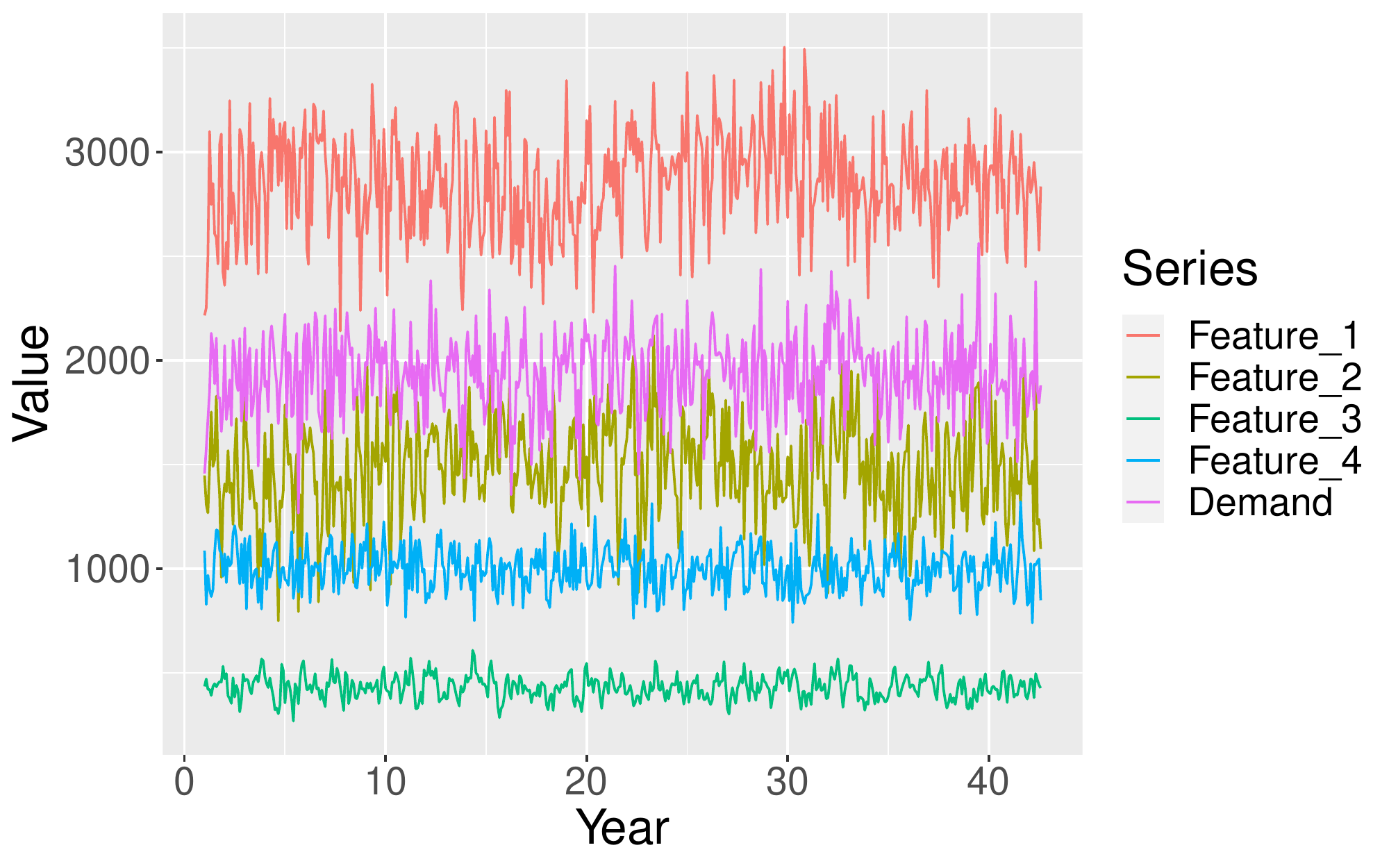}
\end{figure}

\begin{table}[htb]
\caption{Baseline experiment parameters} \label{tab:base_par}
\centering
\resizebox{\linewidth}{!}{
\begin{tabular}{cccccccccc}
\hline
$b_0$ & $b_1$ & $b_2$ & $b_3$ & $b_4$ & $\theta_{1,1}^1$ & $\theta_{12,1}^1$ & $\theta_{1,1}^2$ & $\phi_{1,1}^2$ & $\theta_{12,1}^2$\\
500 & 0.642 & 0.354 & 0.407 & 0.521 & 0.3 & 0.5 & 0.2 & 0.5 & 0.1\\
\noalign{\smallskip} \hline \noalign{\smallskip}
$\theta_{1,1}^3$ & $\phi_{1,1}^3$ & $\phi_{12,1}^3$ & $\phi_{1,1}^4$ & $\phi_{1,2}^4$ & $\theta_{12,1}^4$ & $\theta_{12,2}^4$ & \texttt{rnorm} & \texttt{rlaplace} & \texttt{rt}\\
0.3 & 0.2 & 0.1 & 0.1 & 0.2 & 0.1 & 0.1 & $\mu=0$ & $\mu=0$ & $\mu=0$\\
& & & & & & & $\sigma=100$ & $b=71$ & $\sigma=100$\\
& & & & & & & & & $\nu=5$\\
\hline
\end{tabular}}
\end{table}
As shown in Table \ref{tab:base_par}, We mark that the $t$-distribution is believed to have heavy-tail, and the normal distribution is believed to have light-tail. We use a mix of $t$-distribution, normal distribution and Laplace distribution with random weights to simulate real circumstance where we have no information about the shape of error distribution in prior. The time series can be seen in Figure \ref{fig:base_series}.

\section{Histogram example} \label{app:hist}

\begin{figure}[htb]
\centering
\caption{Profit histogram on the non-parametric method and benchmark method}
\label{fig:npsa_hist}
\includegraphics[width=0.7\textwidth]{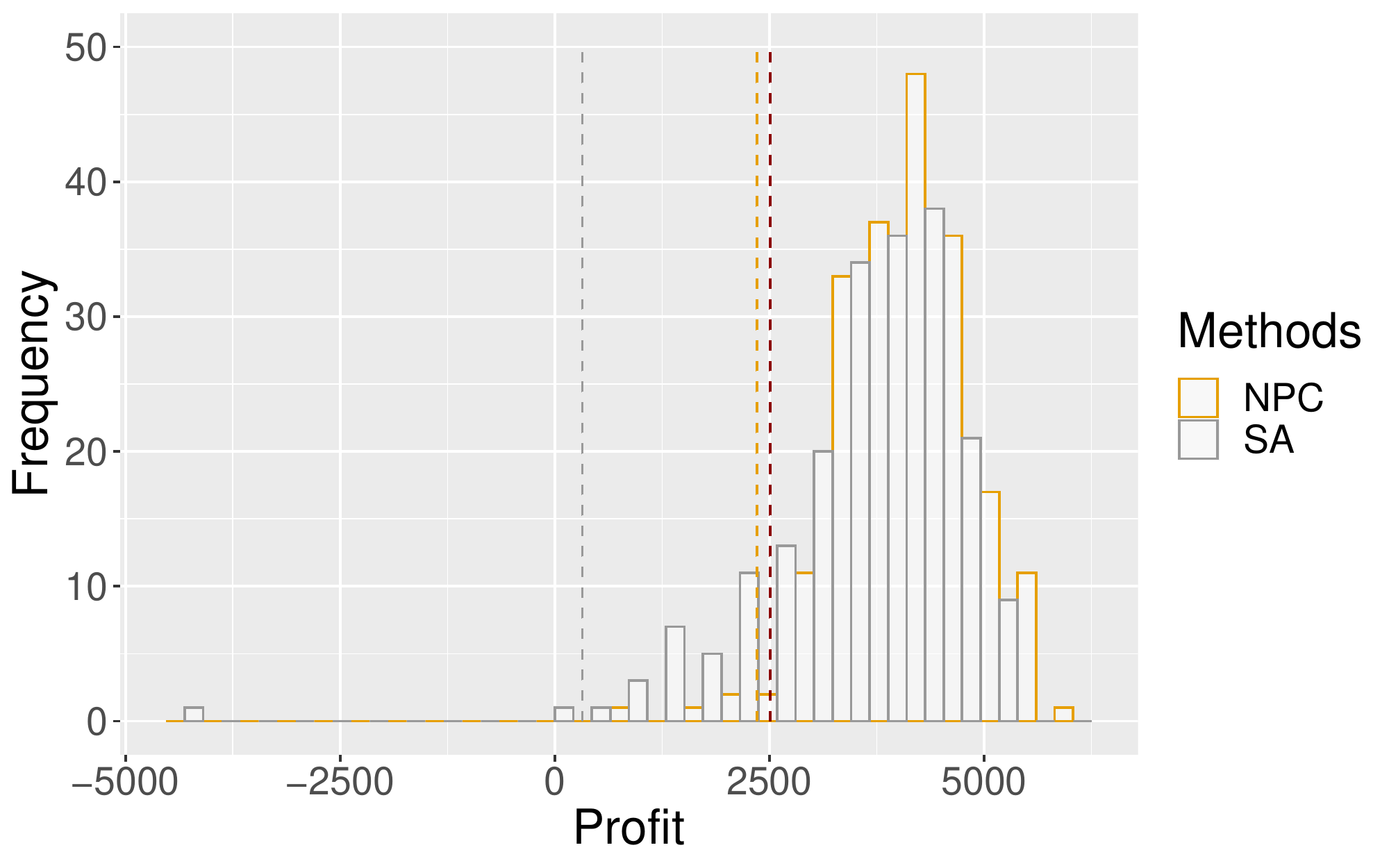}
\end{figure}
In this example, shown in Figure \ref{fig:npsa_hist}, the origin size = 200, iteration number = 100, and $\beta$ = 0.95. Therefore, in the histogram, 100 profit realisations are considered for each methods, and the average profit/loss for the worst 5 cases are marked by dashed lines. We also marked the corresponding performance of the UM method, where all  features and all observations are considered, by a red line dashed line. Here, we have $DL_{SA}=-319.76$, $DL_{NPC}=-2355.24$, $DL_{UM}=-2508.41$, $SL_{SA}=86\%$, $SL_{NPC}=90.5\%$ and $SL_{UM}=88.5\%$.

\newpage
\section{Baseline experiment full results} \label{app:base_results}

\begin{table}[htb]
\caption{Relative $\beta$-$DL$/Relative $SL$ for all choices of parameters under linear profit function when $\beta=0.95$ (or 0.9)} \label{tab:linear_all}
\centering
\resizebox{\linewidth}{!}{
\begin{tabular}{cccccccc}
\hline
\multicolumn{2}{c}{Relative $\beta$-$DL$/$SL$} & \multicolumn{6}{c}{Origin size}\\
\hline
Iteration & Method & 50 & 100 & 150 & 200 & 250 & 300\\
\hline
50 & NPC & 60\%/60\% & 60\%/50\% & 76\%/80\% & 99\%/33\% & 74\%/20\% & 68\%/-\\
& & (55\%/0\%) & (78\%/100\%) & (95\%/25\%) & (89\%/33\%) & (58\%/33\%) & (83\%/0\%)\\
& NF & -39\%/0\% & -1\%/0\% & 11\%/40\% & -13\%/67\% & 8\%/-100\% & -14\%/-\\
& & (-55\%/0\%) & (-5\%/25\%) & (0\%/0\%) & (3\%/-67\%) & (-48\%/-100\%) & (-3\%/-33\%)\\
& SQR & 68\%/20\% & 88\%/25\% & 98\%/60\% & 99\%/99\% & 89\%/0\% & 74\%/-\\
& & (84\%/0\%) & (68\%/75\%) & (99\%/75\%) & (99\%/67\%) & (90\%/67\%) & (95\%/67\%)\\
& PLM & -8\%/40\$ & 7\%/25\% & 7\%/60\% & 99\%/33\% & 57\%/-300\% & 32\%/-\\
& & (-88\%/-60\%) & (72\%/75\%) & (89\%/25\%) & (77\%/100\%) & (37\%/67\%) & (81\%/100\%)\\
\hline
100 & NPC & 26\%/17\% & 54\%/0\% & 87\%/6\% & 99\%/15\% & 76\%/5\% & 86\%/67\%\\
& & (60\%/75\%) & (97\%/63\%) & (93\%/14\%) & (81\%/20\%) & (75\%/0\%) & (91\%/33\%)\\
& NF & -73\%/0\% & 0\%/200\% & 5\%/67\% & -6\%/0\% & -5\%/0\% & 20\%/6\%\\
& & (-19\%/50\%) & (-2\%/13\%) & (-1\%/-29\%) & (-10\%/-100\%) & (-21\%/-67\%) & (5\%/0\%)\\
& SQR & 79\%/-17\% & 96\%/-50\% & 99\%/89\% & 98\%/100\% & 69\%/-200\% & 85\%/67\%\\
& & (79\%/50\%) & (77\%/75\%) & (94\%/71\%) & (99\%/60\%) & (91\%/100\%) & (99\%/50\%)\\
& PLM & -71\%/-33\% & 49\%/-60\% & 80\%/56\% & 97\%/-100\% & 44\%/-20\% & 59\%/89\%\\
& & (4\%/-75\%) & (81\%/75\%) & (79\%/29\%) & (67\%/60\%) & (66\%/50\%) & (88\%/67\%)\\
\hline
150 & NPC & 65\%/43\% & 68\%/50\% & 90\%/14\% & 99\%/25\% & 80\%/56\% & 90\%/56\%\\
& & (77\%/20\%) & (93\%/30\%) & (90\%/8\%) & (83\%/14\%) & (83\%/0\%) & (93\%/44\%)\\
& NF & -25\%/29\% & -11\%/0\% & 6\%/86\% & -1\%/0\% & 19\%/22\% & 12\%/67\%\\
& & (-30\%/30\%) & (-2\%/-20\%) & (-7\%/-25\%) & (-9\%/-85\%) & (-12\%/-10\%) & (19\%/0\%)\\
& SQR & 87\%/-71\% & 92\%/-50\% & 98\%/71\% & 81\%/50\% & 88\%/44\% & 82\%/44\%\\
& & (74\%/60\%) & (83\%/70\%) & (92\%/83\%) & (97\%/71\%) & (97\%/80\%) & (99\%/78\%)\\
& PLM & -45\%/-44\% & 57\%/-50\% & 79\%/43\% & 86\%/-100\% & 68\%/56\% & 64\%/67\%\\
& & (-21\%/-20\%) & (74\%/100\%) & (75\%/75\%) & (69\%/14\%) & (75\%/60\%) & (83\%/100\%)\\
\hline
200 & NPC & 76\%/72\% & 68\%/0\% & 77\%/0\% & 73\%/44\% & 71\%/73\% & 90\%/20\%\\
& & (71\%/40\%) & (91\%/73\%) & (90\%/13\%) & (87\%/36\%) & (84\%/15\%) & (87\%/5\%)\\
& NF & -24\%/29\% & -8\%/-80\% & 0\%/80\% & -8\%/11\% & -5\%/27\% & -14\%/100\%\\
& & (-21\%/20\%) & (-2\%/-13\%) & (-8\%/-7\%) & (-7\%/-18\%) & (4\%/-8\%) & (14\%/16\%)\\
& SQR & 70\%/-43\% & 77\%/-50\% & 75\%/40\% & 78\%/78\% & 76\%/55\% & 72\%/0\%\\
& & (71\%/73\%) & (82\%/67\%) & (92\%/93\%) & (99\%/64\%) & (97\%/92\%) & (99\%/95\%)\\
& PLM & 42\%/-43\% & 38\%/-40\% & 44\%/40\% & 39\%/44\% & 51\%/73\% & 47\%/80\%\\
& & (-14\%/-7\%) & (60\%/93\%) & (76\%/87\%) & (71\%/45\%) & (68\%/77\%) & (88\%/95\%)\\
\hline
\end{tabular}}
\end{table}

\newpage
\begin{table}[htb]
\caption{Relative $\beta$-$DL$/Relative $SL$ for all choices of parameters under nonlinear profit function when $\beta=0.95$ (or 0.9)} \label{tab:nonlinear_all}
\centering
\resizebox{\linewidth}{!}{
\begin{tabular}{cccccccc}
\hline
\multicolumn{2}{c}{Relative $\beta$-$DL$/$SL$} & \multicolumn{6}{c}{Origin size}\\
\hline
Iteration & Method & 50 & 100 & 150 & 200 & 250 & 300\\
\hline
50 & NPC & 85\%/20\% & 99\%/50\% & 56\%/- & 66\%/88\% & 99\%/25\% & 99\%/100\%\\
& & (87\%/100\%) & (78\%/50\%) & (63\%/0\%) & (75\%/50\%) & (95\%/-) & (85\%/90\%)\\
& NF & -17\%/-100\% & -84\%/0\% & 8\%/- & 6\%/13\% & 13\%/0\% & -14\%/11\%\\
& & (-39\%/-100\%) & (-10\%/-50\%) & (11\%/50\%) & (4\%/13\%) & (4\%/-) & (-9\%/20\%)\\
& SQR & 1\%/-40\% & 5\%/0\% & 6\%/- & 23\%/0\% & 43\%/-50\% & 37\%/-11\\
& & (-8\%/-40\%) & (-16\%/-50\%) & (-3\%/-15\%) & (12\%/13\%) & (18\%/-) & (18\%/-20\%)\\
& PLM & -10\%/-100\$ & -2\%/-50\% & 39\%/- & 36\%/88\% & 63\%/100\% & 57\%/88\\
& & (-32\%/100\%) & (-19\%/0\%) & (52\%/100\%) & (20\%/88\%) & (52\%/-) & (72\%/100\%)\\
\hline
100 & NPC & 99\%/60\% & 99\%/100\% & 87\%/40\% & 80\%/58\% & 84\%/27\% & 99\%/57\%\\
& & (99\%/67\%) & (97\%/100\%) & (83\%/80\%) & (92\%/25\%) & (88\%/60\%) & (91\%/14\%)\\
& NF & -43\%/-100\% & -19\%/0\% & 23\%/0\% & 9\%/8\% & -3\%/0\% & -12\%/-14\%\\
& & (-51\%/33\%) & (-31\%/0\%) & (6\%/0\%) & (-7\%/-25\%) & (0\%/-10\%) & (-3\%/42\%)\\
& SQR & -32\%/-50\% & -15\%/-25\% & 0\%/14\% & 43\%/-8\% & 46\%/-36\% & 19\%/-85\%\\
& & (-15\%/-13\%) & (-19\%/-20\%) & (-3\%/-14\%) & (22\%/-37\%) & (33\%/-50\%) & (32\%/-14\%)\\
& PLM & -21\%/-19\% & -23\%/-45\% & 39\%/80\% & 53\%/91\% & 74\%/72\% & 62\%/14\%\\
& & (4\%/33\%) & (81\%/67\%) & (79\%/100\%) & (67\%/63\%) & (66\%/70\%) & (88\%/85\%)\\
\hline
150 & NPC & 85\%/30\% & 99\%/14\% & 86\%/22\% & 99\%/17\% & 99\%/100\% & 88\%/13\%\\
& & (74\%/90\%) & (93\%/50\%) & (88\%/100\%) & (91\%/23\%) & (93\%/83\%) & (67\%/50\%)\\
& NF & -16\%/-100\% & -13\%/40\% & 14\%/11\% & -6\%/11\% & -3\%/12\% & -11\%/-33\%\\
& & (-37\%/100\%) & (-21\%/100\%) & (2\%/-100\%) & (-5\%/-5\%) & (0\%/-33\%) & (-6\%/50\%)\\
& SQR & -13\%/-90\% & -10\%/-18\% & 27\%/-88\% & 33\%/-23\% & 24\%/12\% & 10\%/18\%\\
& & (22\%/-11\%) & (6\%/-46\%) & (2\%/-30\%) & (6\%/-29\%) & (3\%/-21\%) & (5\%/-16\%)\\
& PLM & -13\%/-29\% & -10\%/-16\% & 45\%/77\% & 60\%/88\% & 59\%/25\% & 67\%/-67\%\\
& & (3\%/-11\%) & (52\%/66\%) & (59\%/75\%) & (49\%/76\%) & (43\%/50\%) & (51\%/88\%)\\
\hline
200 & NPC & 76\%/10\% & 56\%/18\% & 78\%/69\% & 96\%/67\% & 75\%/21\% & 93\%/85\%\\
& & (60\%/87\%) & (95\%/75\%) & (77\%/80\%) & (88\%/83\%) & (97\%/83\%) & (97\%/85\%)\\
& NF & -24\%/50\% & -6\%/-33\% & 1\%/7\% & -9\%/-16\% & -4\%/14\% & -12\%/-15\%\\
& & (-34\%/25\%) & (-20\%/-50\%) & (-2\%/-8\%) & (-6\%/15\%) & (-4\%/-16\%) & (-6\%/33\%)\\
& SQR & -8\%/-55\% & 12\%/-23\% & 25\%/-76\% & 16\%/-91\% & 17\%/-21\% & 4\%/-11\%\\
& & (-16\%/-35\%) & (-38\%/11\%) & (9\%/-93\%) & (18\%/-92\%) & (30\%/-30\%) & (-19\%/-88\%)\\
& PLM & -12\%/-16\% & -5\%/-100\% & 46\%/61\% & 42\%/91\% & 52\%/-42\% & -2\%/-15\%\\
& & (-9\%/-20\%) & (55\%/-10\%) & (57\%/75\%) & (30\%/76\%) & (37\%/50\%) & (9\%/66\%)\\
\hline
\end{tabular}}
\end{table}

\newpage
\section{Features for food preparation problem} \label{app:invexp}
\begin{table}[ht!]
\caption{Relevant features to food preparation problem
within food bank} \label{tab:food_features}
\centering
\resizebox{0.5\linewidth}{!}{
\begin{tabular}{cc}
\hline
No. & Feature\\
\hline
1 & UK inflation data (monthly)\\
2 & UK unemployment rate (monthly)\\
3 & UK economics index (weekly)\\
4 & FTSE 100 close price (weekly)\\
5 & Durham birth registered (weekly)\\
6 & Durham death registered (weekly)\\
7 & Durham Covid-19 cases (weekly)\\
8 & Durham crime index (weekly)\\
9 & UK Bank holidays dummies\\
10 & Seasonality dummies\\
\hline
\end{tabular}}
\end{table}

\newpage
\bibliography{newsvendor4}

\end{document}